\documentclass[11pt]{amsart}
\usepackage[utf8]{inputenc}
\usepackage[T1]{fontenc}
\usepackage{amssymb,amsmath,mathtools}
\usepackage{enumitem}
\usepackage{booktabs}
\usepackage[margin=1.2in]{geometry}
\usepackage[colorlinks=true,linkcolor=blue!50!black,citecolor=green!35!black,urlcolor=blue!50!black]{hyperref}
\usepackage{xcolor}

\newtheorem{theorem}{Theorem}[section]
\newtheorem{lemma}[theorem]{Lemma}
\newtheorem{proposition}[theorem]{Proposition}
\newtheorem{corollary}[theorem]{Corollary}
\newtheorem{conjecture}[theorem]{Conjecture}
\newtheorem{question}[theorem]{Question}
\newtheorem{hypothesis}{Hypothesis}

\newtheorem*{theoremA}{Theorem A}
\newtheorem*{theoremF}{Theorem F}
\newtheorem*{theoremS}{Theorem S}
\newtheorem*{theoremW}{Theorem W}
\newtheorem*{theoremP}{Theorem P}
\newtheorem*{theoremM}{Theorem M}
\newtheorem*{theoremO}{Theorem O}
\newtheorem*{conjectureI}{Conjecture I}
\newtheorem*{conjectureII}{Conjecture II}
\theoremstyle{definition}
\newtheorem{definition}[theorem]{Definition}
\newtheorem{example}[theorem]{Example}
\newtheorem{computation}[theorem]{Computation}
\theoremstyle{remark}
\newtheorem{remark}[theorem]{Remark}
\newtheorem{observation}{Observation}[section]

\newcommand{\R}{\mathbb{R}}
\newcommand{\Z}{\mathbb{Z}}
\newcommand{\F}{\mathbb{F}}
\newcommand{\AM}{\mathcal{M}}
\newcommand{\Sep}{\mathcal{S}}
\newcommand{\Gr}{\mathcal{G}}
\newcommand{\fat}{\mathfrak{M}}
\newcommand{\Pl}{\mathfrak{P}}
\newcommand{\V}{\mathcal{V}}
\newcommand{\Cc}{\mathcal{C}}
\newcommand{\Rad}{\operatorname{Rad}}
\newcommand{\supp}{\operatorname{supp}}
\newcommand{\sgn}{\operatorname{sign}}
\newcommand{\one}{\mathbf{1}}
\newcommand{\Sph}{\Sigma}

\newcommand{\Haus}{\mathcal{K}}
\newcommand{\Sp}{\mathsf{S}}
\newcommand{\T}{\mathsf{T}}
\newcommand{\rad}{\mathrel{\dagger}}
\newcommand{\sep}{\mathrel{|}}
\newcommand{\up}{\operatorname{up}}

\title{Spherical separoids and the acyclic MacPhersonian}
\author{Ricardo Strausz}
\address{Instituto de Matem\'aticas, Universidad Nacional Aut\'onoma de M\'exico, Ciudad Universitaria, 04510 Ciudad de M\'exico, Mexico}
\email{dino@math.unam.mx}

\begin{document}

\begin{abstract}
The acyclic MacPhersonian $\AM_n^d$ is the poset of acyclic oriented matroids of rank $d+1$ on $n$ labelled elements, the combinatorial model of the space of configurations of $n$ points spanning $\R^d$ up to affine equivalence, which is the Grassmannian $G(d,n-1)$. We study it inside the larger poset $\Sep_n^d$ of \emph{spherical separoids}, those acyclic separoids whose Radon complex is a sphere of the dimension of the space of affine dependences.

Our main unconditional results are as follows. Over every poset of spherical separoids that contains the realizable oriented matroids there is a tautological sphere fibration, and its Stiefel--Whitney classes show that the mod~$2$ cohomology ring of $G(d,n-1)$ is a ring retract of that of the poset, for all $n$ and $d$. For $d=n-2$ the poset of spherical separoids is homotopy equivalent to $G(n-2,n-1)$, and contains $\AM_n^{n-2}$ as a homotopy equivalent subposet. A fat model of $\AM_n^d$ made of piecewise geodesic antipodal spheres in the round sphere is weakly homotopy equivalent to its order complex. The chirotope of every such sphere is recovered from the degrees of its coordinate projections.

We then analyse a strategy based on the mean curvature flow. Any argument that follows the combinatorial type of a sphere generically along the flow would connect every uniform oriented matroid to a realizable one by mutations. By a theorem of Gu and Knauer this is impossible in every rank $r\ge7$ on $n\ge r+17$ elements, so the flow must leave the oriented matroids. We show that the basin hypothesis of the flow, together with the property that $\|\AM_n^d\|$ is a homology manifold, implies $\|\AM_n^d\|\simeq G(d,n-1)$, without following types along the flow. We conclude with the open problems and conjectures that remain.
\end{abstract}

\maketitle

\setcounter{tocdepth}{1}
\tableofcontents

% =====================================================================
\part*{I. Introduction}
\section{Introduction}\label{sec:intro}

\subsection{The MacPhersonian}
Gelfand and MacPherson \cite{GM92} gave a combinatorial formula for the rational Pontryagin classes of a triangulated manifold, in which oriented matroids play the role of tangent planes. Pursuing this idea, MacPherson \cite{Mac93} introduced \emph{combinatorial differential manifolds} and their structure space, the \emph{MacPhersonian} $\mathrm{MacP}(k,n)$: the poset of oriented matroids of rank $k$ on $[n]$, ordered by weak maps. Taking oriented matroids of subspaces gives a stratification of the real Grassmannian $G(k,n)$ indexed by $\mathrm{MacP}(k,n)$, and hence a map $G(k,n)\to\|\mathrm{MacP}(k,n)\|$ to the order complex. The problem of whether this map is a homotopy equivalence was popularized by Mn\"ev and Ziegler \cite{MZ93}. They also showed that the corresponding statement fails for the Grassmannians of non-realizable oriented matroids.

The difficulty is Mn\"ev's universality theorem \cite{Mnev88} (see \cite{Vakil06} for a refereed account). The strata of the stratification, i.e.\ realization spaces, can have the homotopy type of any primary semialgebraic set, so the map cannot be studied stratum by stratum.

Here is what is known.
\begin{itemize}
\item In rank $2$ the conjecture was proved by Babson \cite{Bab93}, and a homeomorphism was given by Abawonse \cite{Aba23}.
\item Anderson studied the homotopy groups of $\mathrm{MacP}(k,n)$ \cite{And98}. Anderson and Davis \cite{AD02} proved that the mod~$2$ cohomology of $\|\mathrm{MacP}(k,\infty)\|$ contains the polynomial ring on the Stiefel--Whitney classes, using matroid bundles \cite{And99}. Later they identified $\mathrm{MacP}(k,n)$ with the Grassmannian over the sign hyperfield, and showed that $G(k,n)\to\|\mathrm{MacP}(k,n)\|$ is surjective in mod~$2$ cohomology for finite $n$ and split stably \cite{AD19}.
\item Biss \cite{Bis03} published a proof of the general case. Mn\"ev found an error in it, acknowledged in an erratum \cite{Bis09}, and the problem is open again.
\item Dobbins \cite{Dob21} related Grassmannians to spaces of weighted pseudosphere arrangements, and has announced the rank $3$ case \cite{Dob25}.
\item Inductive approaches via single-element extensions face the existence of oriented matroids, even realizable ones, with non-spherical extension spaces \cite{Liu20}.
\end{itemize}

\subsection{Separoids and the acyclic MacPhersonian}
A \emph{separoid} \cite{Str98,ABMOS02} is a set $E$ with a symmetric relation $A\rad B$ between disjoint subsets, the \emph{Radon partitions}, closed under enlarging either side. The complementary pairs $A\sep B$ are the \emph{separations}. The motivating example is a family of convex sets in $\R^d$, where $A\rad B$ means that the convex hulls of $\bigcup A$ and $\bigcup B$ meet. Every separoid arises in this way \cite{ABMOS02,BS06}. Separoids form a universal category \cite{NS06}, they host Radon-, Tverberg- and Erd\H{o}s--Szekeres-type theorems \cite{MPS06,Str05,Str08}, and they characterize the cocircuit graphs of uniform oriented matroids \cite{MS06}. The \emph{Radon complex} of a separoid was introduced in \cite{MS04} to count polytopes. Oriented matroids are special separoids: $A\rad B$ iff some vector has positive part in $A$ and negative part in $B$.

For a configuration $p_1,\dots,p_n$ spanning $\R^d$, the separoid records which subsets have intersecting convex hulls. Its oriented matroid is acyclic of rank $d+1$. The configuration is determined up to affine equivalence by its space of affine dependences $W\subset\one^\perp\subset\R^n$, of dimension $\ell=n-d-1$. Hence the configuration space is the \emph{affine Grassmannian}
\[\Gr_n^d=G(\ell,\one^\perp)\cong G(\ell,n-1)\cong G(d,n-1).\]
Its combinatorial model is the \emph{acyclic MacPhersonian} $\AM_n^d$ of acyclic oriented matroids of rank $d+1$ on $[n]$. In \cite{Str19} the author proposed to study $\AM_n^d$ through the separoids of such configurations, through a "fat" model of spheres in the round sphere, and through the mean curvature flow, and conjectured the following.

\begin{conjectureI}[acyclic MacPherson conjecture]
$\|\AM_n^d\|\simeq G(d,n-1)$ for all $1\le d\le n-2$.
\end{conjectureI}

\noindent By duality, $\AM_n^d$ is isomorphic to the poset of totally cyclic oriented matroids of rank $n-d-1$. The genuinely new range is $3\le d\le n-5$; the first open case is $(n,d)=(8,3)$, where non-realizable oriented matroids first appear.

The Radon partitions of a configuration span an open subset of the sphere $\Sph=\one^\perp\cap S^{n-1}$, the union of the coordinate cells $C_{A,B}$ of sign pattern $(A,B)$ with $A\rad B$. It contains the great sphere $W\cap\Sph$ as a deformation retract. A separoid is \emph{spherical} if this open set, its Radon complex, has the homotopy type of $S^{\ell-1}$. Let $\Sep_n^d$ be the poset of spherical acyclic separoids of dimension $d$ on $[n]$ (Definition~\ref{def:spherical}). It contains $\AM_n^d$.

\begin{conjectureII}[spherical separoids]
$\|\Sep_n^d\|\simeq G(d,n-1)$, and the inclusion $\AM_n^d\hookrightarrow\Sep_n^d$ is a homotopy equivalence.
\end{conjectureII}

\noindent Conjecture~II holds for $d=n-2$ (Theorem~F below). For $n\le5$, its analogue for spherical \emph{Radon} separoids holds by computer, and for $n=4$ it holds for all spherical separoids, in homology when $d=1$ (Section~\ref{sec:computations}). Without sphericity it fails already for $n=4$.

\subsection{Results}
This paper is organised in four parts:
\begin{itemize}
\item the introduction;
\item results proved in full, some of them with the help of a computer (Sections~\ref{sec:separoids}--\ref{sec:computations});
\item conditional results, i.e.\ proved implications between hypotheses on the mean curvature flow and the conjectures (Section~\ref{sec:conditional});
\item open problems and conjectures (Section~\ref{sec:open}).
\end{itemize}
Numerical experiments are collected separately in Appendix~\ref{app:numerics} and are not used in any proof. We use the notation of separoids throughout.

\begin{theoremA}[Section~\ref{sec:extreme}]
For $d\in\{1,n-2\}$ there is a homeomorphism $\|\AM_n^d\|\cong\R P^{n-2}\cong G(d,n-1)$.
\end{theoremA}

\begin{theoremF}[Section~\ref{sec:extreme}]
Let $n\ge3$. For every poset $\mathcal B$ of homotopy-spherical acyclic separoids on $[n]$ with $\AM_n^{n-2}\subseteq\mathcal B$, in particular for $\mathcal B=\Sep_n^{n-2}$, the inclusion $\AM_n^{n-2}\hookrightarrow\mathcal B$ is a homotopy equivalence and $\|\mathcal B\|\simeq G(n-2,n-1)$. Thus Conjecture~II holds for $d=n-2$.
\end{theoremF}

\begin{theoremW}[Section~\ref{sec:W}]
Let $\mathcal B$ be any poset of $\F_2$-spherical acyclic separoids on $[n]$, of sphere dimension $\ell-1=n-d-2$, that contains the realizable oriented matroids; e.g.\ $\mathcal B=\AM_n^d$ or $\Sep_n^d$. There are classes $w_i(\mathcal B)\in H^i(\|\mathcal B\|;\F_2)$, the Stiefel--Whitney classes of the tautological sphere fibration over $\mathcal B$, which pull back to the Stiefel--Whitney classes of the tautological bundle under the realization map $f:G(d,n-1)\to\|\mathcal B\|$, and which satisfy all the relations of the Grassmannian. Consequently
\[f^*:H^*(\|\mathcal B\|;\F_2)\to H^*(G(d,n-1);\F_2)\]
is a split surjection of rings.
\end{theoremW}

\noindent Sphericity is exactly what makes the tautological fibration exist. Acyclicity produces the relations $\bar w_j=0$ for $j>d$: all Radon partitions are mixed, so every Radon complex lives in the $(n-2)$-sphere $\Sph\subset\one^\perp$. As consequences, the mod~$2$ form of Conjecture~II is equivalent to the numerical bound $\sum_ib_i(\|\Sep_n^d\|;\F_2)\le\binom{n-1}{d}$. Moreover, if $\|\AM_n^d\|$ is an $\F_2$-homology manifold of dimension $d(n-1-d)$, then $H^*(\|\AM_n^d\|;\F_2)\cong H^*(G(d,n-1);\F_2)$ (Corollaries~\ref{cor:betti} and~\ref{cor:manifold}).

\begin{theoremS}[Section~\ref{sec:fat}]
Let $\fat_n^d$ be the space of antipodally symmetric piecewise geodesic $(\ell-1)$-spheres in $\Sph$ whose cells meet exactly the coordinate cells of the vectors of an acyclic oriented matroid, with the Hausdorff topology. The separoid map $\fat_n^d\to\AM_n^d$ is a weak homotopy equivalence; in particular $\fat_n^d\simeq_w\|\AM_n^d\|$.
\end{theoremS}

\noindent Under the separoid map, the Quillen fibre of $\AM_n^d\hookrightarrow\Sep_n^d$ over $\Sp$ corresponds to the space of spheres of $\fat_n^d$ contained in the Radon complex of $\Sp$ (Proposition~\ref{prop:fibres-spheres}).

\begin{theoremP}[Section~\ref{sec:degree}]
Let $\Sph'\in\fat_n^d$ be an oriented sphere of type $M$, and let $B\subseteq[n]$ with $|B|=\ell$.
\begin{itemize}
\item If the dual chirotope $\chi^*_M(B)$ is nonzero, the projection $\Sph'\to S^{\ell-1}$, $x\mapsto x_B/|x_B|$, has odd degree, and its sign is $\varepsilon\,\chi^*_M(B)$ for a sign $\varepsilon$ independent of $B$.
\item If $\chi^*_M(B)=0$, the sphere meets the coordinate subspace $\{x_B=0\}$.
\end{itemize}
\end{theoremP}

\noindent Theorem~P relates $\fat_n^d$ to a second fat model, the \emph{Pl\"ucker fat MacPhersonian} $\Pl\subseteq S^{\binom n\ell-1}$ (Section~\ref{sec:degree}). It is the union of the open orthant cells of the Pl\"ucker sphere indexed by the dual chirotopes of acyclic oriented matroids. Its open stars are star-shaped, so the star lemma is free there, and $\Pl/\pm\simeq\|\AM_n^{d,\chi}\|$ for the weak-map order. The Plücker-embedded Grassmannian meets $\Pl$ in its acyclic locus, which is homotopy equivalent to $G(d,n-1)$. We show that the acyclic MacPherson conjecture for the weak-map order is equivalent to the inclusion $G(d,n-1)\hookrightarrow\Pl/\pm$ being a weak equivalence (Theorem~\ref{thm:plucker-equiv}). In words: \emph{replacing the Pl\"ucker relations by their sign shadows, the $3$-term Grassmann--Pl\"ucker conditions, should not change the homotopy type.} The degrees of Theorem~P define a map $\overline\Phi:\fat_n^d\to\Pl/\pm$ that restricts to the Pl\"ucker embedding on great spheres.

\subsection{The flow, and why oriented matroids are not enough}\label{subsec:obstruction}
The strategy of \cite{Str19} realizes every acyclic oriented matroid by an antipodal sphere $\Sph'\in\fat_n^d$ and lets the mean curvature flow of the round sphere deform it. Since the flow is odd, it cannot collapse $\Sph'$ to a point, and one hopes that it converges to a great sphere $W\cap\Sph$, $W\in\Gr_n^d$ (the basin hypothesis~\ref{H:basin}). If the flow could also be followed combinatorially (\ref{H:type}), then Theorem~S would give Conjecture~I (Theorem~\ref{thm:B}).

In numerical experiments at $(8,3)$ (Appendix~\ref{app:numerics}), starting from certified non-realizable oriented matroids, the combinatorial type changes only by \emph{mutations}, i.e.\ by reversing the sign of a single basis. This suggests the most natural mechanism for following the flow: a generic trajectory crosses only codimension-one walls, which are mutations. Our first conditional result shows that this mechanism cannot work in general.

\begin{theoremO}[Section~\ref{sec:conditional}]
Assume the basin hypothesis~\ref{H:basin}, the degree hypothesis~\ref{H:degree} (the degree pattern of Theorem~P remains the chirotope of an acyclic oriented matroid along the flow) and genericity of trajectories. Then every uniform oriented matroid of rank $d+1$ on $[n]$ is joined to a realizable one by a sequence of mutations through uniform oriented matroids. In particular, the mutation graph of uniform oriented matroids of rank $d+1$ on $[n]$ is connected.
\end{theoremO}

Connectivity of the mutation graph is a conjecture of Cordovil and Las Vergnas. Gu and Knauer \cite{GK26} have shown that it is false: for every rank $r\ge7$ and every $n\ge r+17$ the mutation graph of uniform oriented matroids of rank $r$ on $n$ elements is disconnected. In the Pl\"ucker model this says that $\Pl$ contains \emph{islands}: uniform cells of non-realizable oriented matroids that are joined to the realizable cells only through cells of codimension at least $2$ (Proposition~\ref{prop:islands}). Since a generic path crosses only walls of codimension one, a generic deformation cannot reach the Grassmannian from an island while staying in $\Pl$. Hence:

\begin{corollary}[Corollary~\ref{cor:GK}]\label{cor:intro-GK}
For $d\ge6$ and $n\ge d+18$, Hypotheses~\ref{H:basin}, \ref{H:degree} and genericity cannot hold simultaneously.
\end{corollary}

This is the central obstruction of the paper, and it shapes everything that follows. Any proof of Conjecture~I through the flow must let the evolving sphere leave the oriented matroids. Then the combinatorial shadow of the sphere is no longer an oriented matroid. Our experiments (Appendix~\ref{app:numerics}) show the flow doing exactly this at $(8,3)$: the degree pattern violates a Grassmann--Pl\"ucker relation, and the sphere's separoid is not even a Radon separoid. In every check, however, the separoid remains \emph{spherical}. This is why spherical separoids, and not oriented matroids, are the central objects here. It is also why we look for conditional theorems that require the flow only to converge, without following types:
\begin{itemize}
\item Theorem~M (Section~\ref{sec:conditional}): \ref{H:basin} together with the property that $\|\AM_n^d\|$ is a $\Z$-homology manifold of dimension $d(n-1-d)$ implies Conjecture~I. The flow supplies a retraction onto $\Gr_n^d$, and Poincar\'e duality supplies the rest.
\item Theorem~\ref{thm:Eflat} and Proposition~\ref{prop:factor}(c): weaker substitutes for \ref{H:degree}, in terms of contractible fibres or of spherical separoids.
\item Proposition~\ref{prop:factor}(a),(b): inside the space of spherical separoids, the flow factors the inclusion of $\AM_n^d$ through the Grassmannian.
\end{itemize}
Section~\ref{sec:conditional} contains a complete table of the proved implications between the hypotheses.

\subsection{Conventions}
Statements proved with the help of a computer are marked as such, and the computations are described in Section~\ref{sec:computations}; the code is available from the author. Section~\ref{sec:conditional} and Appendix~\ref{app:numerics} contain no unconditional claims beyond those explicitly labelled.

% =====================================================================
\part*{II. Proved results}
% =====================================================================
\section{Separoids, oriented matroids and Radon complexes}\label{sec:separoids}

\subsection{Signed sets}
Let $E=[n]$. A \emph{signed set} is a pair $Z=(Z^+,Z^-)$ of disjoint subsets of $E$; we identify it with its sign vector in $\{+,-,0\}^E$, write $\supp Z=Z^+\cup Z^-$ and $-Z=(Z^-,Z^+)$. The \emph{conformal order} is $Z\le Z'$ iff $Z^+\subseteq Z'^+$ and $Z^-\subseteq Z'^-$. The composition $Z\circ Z'$ agrees with $Z$ on $\supp Z$ and with $Z'$ elsewhere. Let
\[P=P_n=\{Z:\ Z^+\ne\emptyset\ne Z^-\}\]
be the poset of \emph{mixed} signed sets. It is an up-set of the poset of nonzero signed sets, and $Z\mapsto-Z$ is an order-preserving involution of $P$ such that $Z$ and $-Z$ have neither a common upper bound nor a common lower bound.

\subsection{Separoids}
\begin{definition}[\cite{Str98,ABMOS02}]\label{def:separoid}
An \emph{acyclic separoid} on $E$ is a relation $\rad$ on pairs of disjoint subsets of $E$ such that:
\begin{enumerate}[label=(\roman*)]
\item $A\rad B$ implies $B\rad A$;
\item $A\rad B$ and $B\subseteq B'\subseteq E\setminus A$ imply $A\rad B'$;
\item $A\rad B$ implies $A\ne\emptyset\ne B$.
\end{enumerate}
The pairs with $A\rad B$ are the \emph{Radon partitions}. A disjoint pair with $A,B\neq\emptyset$ that is not a Radon partition is a \emph{separation}, written $A\sep B$.
\end{definition}

\noindent Condition (iii) is acyclicity. Without it one obtains general separoids, which play no role here. By (i)--(iii), an acyclic separoid $\Sp$ is the same thing as a symmetric up-set
\[\Rad(\Sp)=\{(A,B):A\rad B\}\subseteq P,\]
and we use both descriptions. Its minimal elements are the \emph{minimal Radon partitions}.
\begin{itemize}
\item The \emph{dimension} of $\Sp$ is the least $d$ such that every $(d+2)$-subset $X\subseteq E$ admits a partition $X=A\sqcup B$ with $A\rad B$.
\item $\Sp$ is a \emph{Radon separoid} if two minimal Radon partitions with the same support coincide up to sign.
\item Separoids are ordered by
\[\Sp\le\T\iff\Rad(\T)\subseteq\Rad(\Sp),\]
i.e.\ every Radon partition of $\T$ is one of $\Sp$; the larger separoid is the more general one.
\end{itemize}

\begin{example}\label{ex:convex}
For points $p_1,\dots,p_n\in\R^d$, the relation $A\rad B$ iff $\mathrm{conv}\,p_A\cap\mathrm{conv}\,p_B\neq\emptyset$ is an acyclic separoid. Its separations are the pairs strictly separated by a hyperplane. By Radon's theorem its dimension is at most $d$, and it is exactly $d$ when the points span $\R^d$ affinely. More generally, any family of convex sets defines a separoid in this way, and every separoid arises from such a family \cite{ABMOS02,BS06}.
\end{example}

\subsection{Oriented matroids as separoids}
We use \cite{BLSWZ} for oriented matroids. For an oriented matroid $M$ on $E$ let $\Cc(M)$ be its signed circuits and $\V(M)$ its nonzero vectors. Vectors are the covectors of the dual $M^*$, and $\V(M)\cup\{0\}$ is closed under composition \cite[\S3.7]{BLSWZ}. The \emph{separoid of $M$} is
\[\Rad(M)=\up\V(M)=\{Z:Z\ge X\text{ for some }X\in\V(M)\}.\]
Its minimal elements are the circuits, since every vector dominates a circuit. So $M$ is determined by its separoid. $M$ is \emph{acyclic} if it has no nonzero vector with empty negative part; then $\Rad(M)\subseteq P$ and $\Rad(M)$ is an acyclic separoid.

\begin{definition}\label{def:AM}
$\AM_n^d$ is the poset of acyclic oriented matroids of rank $d+1$ on $[n]$, ordered as separoids: $M\le M'$ iff $\Rad(M')\subseteq\Rad(M)$, i.e.\ iff every circuit of $M'$ is a Radon partition of $M$.
\end{definition}

\noindent Acyclic oriented matroids have no loops. Their separoids have dimension $d$, because every $(d+2)$-set contains the support of a circuit, and some $(d+1)$-set (a basis) does not. They are Radon separoids, because a matroid circuit has a unique signing up to sign. Other orders on $\AM_n^d$ appear in the literature, namely weak maps and the vector order. They are compared in Section~\ref{sec:conditional}. None of our unconditional results depends on them.

\begin{lemma}\label{lem:ymax}
Let $M$ be an oriented matroid and $Z\in\Rad(M)$. Then $\{Y\in\V(M):Y\le Z\}$ is nonempty and has a maximum $y_M(Z)$. The map $y_M:\Rad(M)\to\V(M)$ is order-preserving, odd, satisfies $y_M(Z)\le Z$, and is the identity on $\V(M)$.
\end{lemma}
\begin{proof}
Nonemptiness is the definition of $\Rad(M)$. If $Y,Y'\le Z$ then $Y\circ Y'\le Z$ and $Y\circ Y'\in\V(M)$. So the composition of all vectors below $Z$ is a vector below $Z$, and it is the maximum. The other properties are immediate.
\end{proof}

\subsection{The affine Grassmannian}
Let $p_1,\dots,p_n$ affinely span $\R^d$ and let $W=\{\lambda\in\R^n:\sum\lambda_ip_i=0,\ \sum\lambda_i=0\}$. It is an $\ell$-dimensional subspace of $\one^\perp$, where
\[\ell=n-d-1,\]
and every such subspace arises. Two configurations are affinely equivalent iff they have the same $W$. The vectors of the oriented matroid $M_W$ of the configuration are the sign vectors of $W\setminus0$, and $M_W\in\AM_n^d$. By Example~\ref{ex:convex}, $\Rad(M_W)$ is the separoid of the configuration. We put
\[\Gr_n^d=G(\ell,\one^\perp)\cong G(\ell,n-1)\cong G(d,n-1),\]
the last isomorphism by orthogonal complements in $\one^\perp$. The \emph{realizable} oriented matroids $M_W$ form a subposet $\AM_n^{d,\mathrm{re}}\subseteq\AM_n^d$.

\subsection{The sphere and its cells}
Let $\Sph=\one^\perp\cap S^{n-1}$, a round sphere of dimension $n-2$. Every $x\in\Sph$ has coordinates of both signs, so
\[\Sph=\bigsqcup_{Z\in P}C_Z,\qquad C_Z=\{x\in\Sph:\sgn x=Z\}.\]
For $Z\in P$ let $\chi_Z\in\{0,\pm1\}^n$ be $Z$ read as a real vector.

\begin{lemma}\label{lem:cells}
\begin{enumerate}[label=(\alph*)]
\item $\overline{C_Z}=\bigsqcup_{Z'\le Z}C_{Z'}$, and $\langle x,\chi_Z\rangle=\sum_i|x_i|>0$ on $\overline{C_Z}$.
\item The gnomonic map $\gamma(x)=2x/\sum_i|x_i|$ maps $\overline{C_Z}$ homeomorphically onto the polytope $Q_Z=\{y\in\one^\perp:\sgn y\le Z,\ \sum_{i\in Z^+}y_i=1\}\cong\Delta^{|Z^+|-1}\times\Delta^{|Z^-|-1}$, and $C_Z$ onto its relative interior.
\item On $\overline{C_Z}$, $\gamma$ is the central projection onto the hyperplane $\{\langle y,\chi_Z\rangle=2\}$. It maps arcs of great circles in $\overline{C_Z}$ to segments, so $\overline{C_Z}$ is geodesically convex, and $\dim C_Z=|\supp Z|-2$.
\end{enumerate}
In particular $\{C_Z\}_{Z\in P}$ is a regular cell decomposition of $\Sph$ with face poset $P$, invariant under $x\mapsto-x$, which maps $C_Z$ to $C_{-Z}$.
\end{lemma}
\begin{proof}
For $x\in\overline{C_Z}$ we have $\langle x,\chi_Z\rangle=\sum|x_i|$, which is positive since $x\neq0$. The point $y=\gamma(x)$ satisfies $\sum|y_i|=2$ and $\sum y_i=0$, hence $\sum_{Z^+}y_i=1=-\sum_{Z^-}y_i$. Conversely, normalizing a point of $Q_Z$ gives a point of $\overline{C_Z}$. A great circle is the intersection of $\Sph$ with a $2$-plane through $0$, and its image under central projection is a line.
\end{proof}

\subsection{Radon complexes and sphericity}

\begin{lemma}\label{lem:upset}
Let $X$ be a regular CW complex with face poset $Q$, and $U\subseteq Q$ an up-set. The union $X_U$ of the open cells of $U$ is open, and it deformation retracts onto the subcomplex $\|U\|$ of the barycentric subdivision $\|Q\|\cong X$. The complement $X\setminus X_U$ is the subcomplex with face poset $Q\setminus U$. If a cellular involution preserves $U$, the retraction can be chosen equivariant.
\end{lemma}
\begin{proof}
A point of $\|Q\|$ lies in the open simplex of a unique chain $q_0<\dots<q_m$, which lies in the open cell of $q_m$. It belongs to $X_U$ iff $q_m\in U$, iff some vertex of the chain lies in $U$, since $U$ is an up-set. Every simplex meeting $X_U$ is thus the join of a face with vertices in $Q\setminus U$ and a nonempty face with vertices in $U$. Killing linearly the barycentric coordinates of the vertices outside $U$ is a deformation retraction of $X_U$ onto $\|U\|$, which commutes with simplicial automorphisms preserving $U$.
\end{proof}

\begin{definition}\label{def:spherical}
Let $\Sp$ be an acyclic separoid on $[n]$. Its \emph{Radon complex} and \emph{separation complex} are
\[R(\Sp)=\bigcup_{(A,B)\in\Rad(\Sp)}C_{A,B}\subseteq\Sph,\qquad L(\Sp)=\Sph\setminus R(\Sp)=\bigcup_{A\sep B}C_{A,B}.\]
$R(\Sp)$ is open and $L(\Sp)$ is a closed subcomplex. By Lemma~\ref{lem:upset}, $R(\Sp)\simeq\|\Rad(\Sp)\|$ equivariantly. For $\ell\ge1$, $\Sp$ is \emph{$(\ell-1)$-spherical} if $R(\Sp)\simeq S^{\ell-1}$, and \emph{$\F$-spherical} if $R(\Sp)$ has the $\F$-homology of $S^{\ell-1}$. We write
\[\Sep_n^d\qquad(\text{resp. }\Sep_n^{d,\F})\]
for the poset of acyclic separoids of dimension $d$ on $[n]$ that are (resp.\ $\F$-) spherical with $\ell=n-d-1$.
\end{definition}

\begin{remark}\label{rem:radon-versions}
In \cite{MS04,Str19} the Radon complex is a subcomplex of the boundary $O_n$ of the $n$-dimensional cross-polytope. Since the open cells of $O_n$ are indexed by nonzero signed sets and $\Rad(\Sp)\subseteq P$ is an up-set there too, the union of the corresponding \emph{open} cells of $O_n$ is also homotopy equivalent to $\|\Rad(\Sp)\|$ (Lemma~\ref{lem:upset}). So the two notions of sphericity agree. The union of the \emph{closed} cells gives a different, wrong notion. For four points in the plane with $p_1=p_2$ and $p_3,p_4$ in general position, $\Rad=\up\{\pm(\{1\},\{2\})\}$, and the closed version is homotopy equivalent to $S^2$ instead of $S^0$.
\end{remark}

\subsection{The separoid of a subset of the sphere}
For $X\subseteq\Sph$ let $\Sp(X)$ be the acyclic separoid with $\Rad(\Sp(X))=\up\{\sgn x:x\in X\}$. In separoid terms, $A\rad B$ iff some $x\in X$ is positive only in $A$ and negative only in $B$. The following properties are immediate.
\begin{enumerate}[label=(\alph*)]
\item $X\subseteq R(\Sp(X))$.
\item For an acyclic separoid $\T$, $\Sp(X)\ge\T$ iff $X\subseteq R(\T)$.
\item $\Sp(W\cap\Sph)=\Rad(M_W)$ for $W\in\Gr_n^d$.
\end{enumerate}
Let $\Haus(\Sph)$ be the space of nonempty compact subsets of $\Sph$ with the Hausdorff metric, and endow every poset with the Alexandrov topology whose open sets are the up-sets.

\begin{lemma}[semicontinuity]\label{lem:semicont}
For every compact $X\subseteq\Sph$ there is a neighbourhood $\mathcal U$ of $X$ in $\Haus(\Sph)$ such that $\Sp(X')\ge\Sp(X)$ for all $X'\in\mathcal U$. Hence $X\mapsto\Sp(X)$ is continuous from $\Haus(\Sph)$ to the finite space of acyclic separoids on $[n]$. Likewise $W\mapsto M_W$ is continuous from $\Gr_n^d$ to $\AM_n^d$, and $(W,v)\mapsto(M_W,\sgn v)$ is continuous on $\{(W,v):v\in W\cap\Sph\}$.
\end{lemma}
\begin{proof}
Otherwise there are $X_j\to X$ and $x_j\in X_j$ with $\sgn x_j\notin\Rad(\Sp(X))$. Passing to a subsequence, $\sgn x_j=Z$ is constant and $x_j\to x\in X$. Then $\sgn x\le Z$, so $Z\in\Rad(\Sp(X))$, a contradiction. The same argument applies to unit vectors of $W_j\to W$, and to pairs $(W_j,v_j)\to(W,v)$, since $\sgn v_j\ge\sgn v$ for $j$ large.
\end{proof}

% =====================================================================
\section{Topological tools}\label{sec:tools}

\subsection{Finite spaces}
A finite poset $X$ is a finite $T_0$ space with the up-sets as open sets. The minimal open set of $x$ is $U_x=X_{\ge x}$.

\begin{lemma}\label{lem:finite}
\begin{enumerate}[label=(\alph*)]
\item There is a natural weak homotopy equivalence $\mu_X:\|X\|\to X$, sending a point to the largest vertex of its support \cite[Thm.~2]{McCord}.
\item If $f,g:X\to Y$ are order-preserving and $f\le g$ pointwise, then $f\simeq g$ \cite{Stong}, \cite[\S1.2]{Barmak}. Hence $U_x$, which has a minimum, is contractible.
\item If $f:T\to X$ is continuous and every $f^{-1}(U_x)$ is weakly contractible, then $f$ is a weak homotopy equivalence \cite[Thm.~6]{McCord}. More generally it suffices that each restriction $f^{-1}(U_x)\to U_x$ be a weak equivalence.
\end{enumerate}
\end{lemma}

\noindent Since $\mu_X$ is a weak equivalence, every continuous map $f:T\to X$ from a CW complex lifts, uniquely up to homotopy, to $g:T\to\|X\|$ with $\mu_Xg\simeq f$. We identify $H^*(X)$ with $H^*(\|X\|)$ via $\mu_X^*$.

\begin{lemma}[coverings]\label{lem:cover}
Let $\iota$ be an order-preserving involution of a finite poset $X$ such that, for every $x$, the elements $x$ and $\iota x$ have neither a common upper bound nor a common lower bound.
\begin{enumerate}[label=(\alph*)]
\item The quotient space $X/\iota$ is the poset with $[x]\le[y]$ iff $y\ge x$ or $y\ge\iota x$, and $q:X\to X/\iota$ is a $2$-sheeted covering map.
\item $\|X\|\to\|X/\iota\|$ is a $2$-sheeted simplicial covering, and $\|X/\iota\|\cong\|X\|/\iota$.
\item The class $w_X\in H^1(X/\iota;\F_2)$ of $q$ corresponds under $\mu^*$ to the class of the covering in (b).
\end{enumerate}
\end{lemma}
\begin{proof}
(a) The smallest $\iota$-invariant open set containing $x$ is $U_x\cup U_{\iota x}$, so the quotient topology is the Alexandrov topology of the stated relation.
\begin{itemize}
\item \emph{Transitivity.} If $[x]\le[y]\le[z]$, say $y\ge x$ or $y\ge\iota x$, and $z\ge y$ or $z\ge\iota y$, then $z\ge x$ or $z\ge \iota x$, because $\iota$ is order-preserving.
\item \emph{Antisymmetry.} If $[x]\le[y]\le[x]$, replace $y$ by $\iota y$ if necessary so that $y\ge x$. Then $x\ge y$, so $x=y$, or $x\ge\iota y$, in which case $y$ is a common upper bound of $y$ and $\iota y$.
\item \emph{Covering.} $q^{-1}(U_{[x]})=U_x\sqcup U_{\iota x}$, and the union is disjoint by the hypothesis on upper bounds. The restriction $q:U_x\to U_{[x]}$ is bijective, since $y,\iota y\in U_x$ would have the common lower bound $x$. It preserves and reflects the order: if $y,y'\in U_x$ and $y'\ge\iota y$, then $y'\ge x$ and $y'\ge\iota x$, which is impossible.
\end{itemize}

(b) A chain $[x_0]<\dots<[x_m]$ has exactly two lifts. Choose a lift $x_m$ of the top; going down, exactly one of $x_i,\iota x_i$ lies below the element already chosen.

(c) $\mu$ is equivariant. So $\|X\|\to\mu^*X$, $u\mapsto([u],\mu(u))$, is a map of $2$-sheeted coverings of $\|X/\iota\|$ that is bijective on fibres, hence an isomorphism. Double coverings of locally path connected, semilocally simply connected spaces (finite spaces are such, by Lemma~\ref{lem:finite}(b)) are classified, naturally under pull-back, by $H^1(-;\F_2)=\operatorname{Hom}(\pi_1(-),\Z/2)$.
\end{proof}

\noindent All involutions below are induced by $Z\mapsto-Z$, and the hypothesis of Lemma~\ref{lem:cover} holds because it holds in $P$. An equivariant continuous map between spaces with free involutions whose quotient maps are coverings exhibits the source covering as the pull-back of the target covering; hence it pulls the class $w$ back to $w$.

\begin{lemma}[fibre lemma]\label{lem:fibre}
Let $\varphi:A\to Q$ be an order-preserving map of finite posets.
\begin{enumerate}[label=(\alph*)]
\item If every fibre $\varphi^{-1}(Q_{\ge q})$ is contractible, then $\|\varphi\|$ is a homotopy equivalence \cite[Prop.~1.6]{Quillen}.
\item If every fibre $\varphi^{-1}(Q_{\ge q})$ is $\F$-acyclic, then $\|\varphi\|$ induces isomorphisms in $\F$-homology.
\end{enumerate}
\end{lemma}
\begin{proof}[Proof of (b)]
Let $C=A\sqcup Q$, with the orders of $A$ and $Q$ and with $q<a$ iff $q\le\varphi(a)$. The map $r:C\to Q$ given by $r|_A=\varphi$ and $r|_Q=\mathrm{id}$ is order-preserving and $r\le\mathrm{id}$. So $\|Q\|\subset\|C\|$ is a homotopy equivalence with inverse $\|r\|$, which restricts to $\|\varphi\|$ on $\|A\|$. Now add the points of $Q$ to $A$ one at a time, along a linear extension of $Q$. When $q$ is added, its link in the poset built so far is $\|Q_{<q}\|*\|\varphi^{-1}(Q_{\ge q})\|$, which is $\F$-acyclic. A cone glued along an $\F$-acyclic subcomplex does not change $\F$-homology (Mayer--Vietoris). Hence $\|A\|\subset\|C\|$, and therefore $\|\varphi\|$, is an $\F$-homology isomorphism.
\end{proof}

\subsection{Free involutions}
\begin{lemma}\label{lem:gysin}
Let $X$ be a compact space with a free involution, which is either a finite simplicial complex on which the involution acts simplicially and $X\to\bar X=X/\iota$ is a covering, or a topological manifold. Let $w\in H^1(\bar X;\F_2)$ be the class of the double cover $X\to\bar X$.
\begin{enumerate}[label=(\alph*)]
\item If $X$ has the $\F_2$-homology of $S^m$, $m\ge0$, then $H^*(\bar X;\F_2)=\F_2[w]/(w^{m+1})$.
\item $X$ is not $\F_2$-acyclic, and $\chi(X)$ is even.
\item If $X'$ is another such space with the $\F_2$-homology of $S^m$, every equivariant map $g:X\to X'$ induces an isomorphism $H^*(\bar X';\F_2)\to H^*(\bar X;\F_2)$, and for $m\ge1$ an isomorphism $H^m(X';\F_2)\to H^m(X;\F_2)$; in particular, if $X\simeq X'\simeq S^m$, then $g$ has odd degree.
\item If $X$ has the $\F_2$-homology of $S^m$, $m\ge1$, every equivariant map from $X$ to a space with free involution homeomorphic to $S^m$ is surjective.
\end{enumerate}
\end{lemma}
\begin{proof}
The covering is the sphere bundle of the line bundle $X\times_\iota\R$, whose mod~$2$ Euler class is $w$. Its Gysin sequence with $\F_2$ coefficients \cite[Thm.~12.2]{MilnorStasheff} is
\[\cdots\to H^{i-1}(\bar X)\xrightarrow{\cdot w}H^i(\bar X)\xrightarrow{\pi^*}H^i(X)\xrightarrow{\pi_!}H^i(\bar X)\xrightarrow{\cdot w}H^{i+1}(\bar X)\to\cdots.\]
If $H^i(X)=H^{i+1}(X)=0$, then $\cdot w:H^i(\bar X)\to H^{i+1}(\bar X)$ is an isomorphism. Since $\bar X$ is finite dimensional, (a) gives $H^i(\bar X)=0$ for $i>m$.
\begin{itemize}
\item \emph{$m\ge1$.} $X$ and $\bar X$ are connected and $\pi^*$ is an isomorphism in degree $0$, so $\pi_!=0$ in degree $0$ and $\cdot w:H^i(\bar X)\to H^{i+1}(\bar X)$ is an isomorphism for $0\le i\le m-2$. In degree $m$ we obtain the exact sequence $0\to H^{m-1}(\bar X)\to H^m(\bar X)\to\F_2\to H^m(\bar X)\to0$, with $H^{m-1}(\bar X)=\F_2w^{m-1}$. This forces $H^m(\bar X)=\F_2w^m$.
\item \emph{$m=0$.} We get $0\to H^0(\bar X)\to\F_2^2\to H^0(\bar X)\to0$, so $H^0(\bar X)=\F_2$.
\item \emph{Part (b).} If $X$ were acyclic, we would get $0\to H^0(\bar X)\to\F_2\to H^0(\bar X)\to0$, which is impossible. The simplices of $X$ come in pairs, so $\chi(X)=2\chi(\bar X)$.
\end{itemize}
(c) $\bar g^*w=w$, so $\bar g^*$ is an isomorphism of the rings in (a). The Gysin sequences are natural. In degree $m$, $\pi^*=0$ and $\pi_!$ is an isomorphism, and $\pi_!\circ g^*=\bar g^*\circ\pi_!$.

(d) If the map missed a pair $\pm p$, the induced map would factor through $(S^m\setminus\{\pm p\})/\pm$, a noncompact connected $m$-manifold with vanishing $H^m(\,;\F_2)$. That would force $w^m=0$ in $H^m(\bar X)$, contradicting (a).
\end{proof}

\subsection{Leray--Hirsch over a finite space}
\begin{lemma}\label{lem:LH}
Let $p:E\to B$ be continuous with $B$ a finite $T_0$ space, $\F$ a field, and $c_0,\dots,c_m\in H^*(E;\F)$ homogeneous classes whose restrictions form a basis of $H^*(p^{-1}(U_b);\F)$ for every $b\in B$. Then $H^*(E;\F)$ is a free $H^*(B;\F)$-module with basis $c_0,\dots,c_m$, via $a\cdot c=p^*(a)\smile c$.
\end{lemma}
\begin{proof}
Let $\Phi_U:H^*(U)\otimes\langle c_i\rangle\to H^*(p^{-1}U)$, $a\otimes c_i\mapsto p^*a\smile c_i|$. We show by induction on $|U|$ that $\Phi_U$ is an isomorphism for every open $U\subseteq B$. Let $b$ be minimal in $U\ne\emptyset$. Then $U_b$, $U'=U\setminus\{b\}$ and $U_b\cap U'=U_b\setminus\{b\}$ are open.
\begin{itemize}
\item $\Phi_{U_b}$ is an isomorphism by hypothesis, since $H^*(U_b)=\F$ (Lemma~\ref{lem:finite}(b)).
\item $\Phi_{U'}$ and $\Phi_{U_b\setminus b}$ are isomorphisms by induction.
\end{itemize}
The Mayer--Vietoris sequences of the open covers $\{U_b,U'\}$ of $U$ and $\{p^{-1}U_b,p^{-1}U'\}$ of $p^{-1}U$ are compatible with $\Phi$, since cup product with the $c_i$ commutes with restriction and with the coboundary. Tensoring over a field preserves exactness, and the five lemma concludes, as in the classical case \cite[Thm.~4D.1]{Hatcher}.
\end{proof}

\begin{lemma}[Segre coefficients]\label{lem:segre}
Let $A$ be a commutative $\F_2$-algebra, $w_1,\dots,w_\ell\in A$, and $R=A[x]/(x^\ell+w_1x^{\ell-1}+\dots+w_\ell)$, a free $A$-module with basis $1,\dots,x^{\ell-1}$. Let $\pi_!:R\to A$ be the coefficient of $x^{\ell-1}$, and define $\bar w_j$ by $\sum_j\bar w_j=(1+w_1+\dots+w_\ell)^{-1}$. Then $\pi_!(x^{\ell-1+j})=\bar w_j$ for $j\ge0$. In particular, if $x^N=0$ in $R$, then $\bar w_j=0$ for all $j\ge N-\ell+1$.
\end{lemma}
\begin{proof}
Put $s_i=\pi_!(x^i)$, so $s_i=0$ for $i<\ell-1$ and $s_{\ell-1}=1$. Applying the $A$-linear map $\pi_!$ to $x^m(x^\ell+\sum_iw_ix^{\ell-i})=0$ gives $\sum_{i=0}^\ell w_is_{m+\ell-i}=0$ (with $w_0=1$) for $m\ge0$. With $r_j=s_{\ell-1+j}$ this reads $\sum_iw_ir_{j-i}=0$ for $j\ge1$, with $r_0=1$ and $r_{<0}=0$. This is the recursion for $(\sum_iw_i)^{-1}$. If $x^N=0$ and $\ell-1+j\ge N$, then $x^{\ell-1+j}=0$.
\end{proof}
% =====================================================================
\section{Spherical separoids and the extreme cases}\label{sec:extreme}

\subsection{Oriented matroids are spherical}
\begin{proposition}\label{prop:OM-spherical}
Let $M\in\AM_n^d$. Then $R(M)\simeq\|\Rad(M)\|\simeq\|\V(M)\|\cong S^{\ell-1}$. In particular $\AM_n^d\subseteq\Sep_n^d$.
\end{proposition}
\begin{proof}
For the inclusion $\V(M)\hookrightarrow\Rad(M)$, the fibre over $Z$ is $\{Y\in\V(M):Y\le Z\}$. It has the maximum $y_M(Z)$ (Lemma~\ref{lem:ymax}), hence is contractible. By Lemma~\ref{lem:fibre}(a), applied to the opposite posets, $\|\Rad(M)\|\simeq\|\V(M)\|$. The poset $\V(M)\cup\{0\}$ is the poset of covectors of the rank $\ell$ oriented matroid $M^*$. By the topological representation theorem \cite{FL78}, \cite[Ch.~4--5]{BLSWZ} it is the face poset of a regular cell decomposition of $S^{\ell-1}$, so $\|\V(M)\|\cong S^{\ell-1}$. Finally $R(M)\simeq\|\Rad(M)\|$ by Lemma~\ref{lem:upset}.
\end{proof}

\noindent Sphericity is not automatic: without it, the poset of all acyclic Radon separoids of dimension $d$ on $[4]$ has $H_2$ of rank $9$ for $d=1$ and $24$ for $d=2$ (Section~\ref{sec:computations}).

\subsection{Duality between Radon partitions and separations}
\begin{proposition}\label{prop:alexander}
Let $\Sp$ be an acyclic separoid on $[n]$ with $1\le d\le n-2$ and $\ell=n-d-1$. Then $\Sp$ is $\F$-spherical (with sphere dimension $\ell-1$) iff its separation complex $L(\Sp)$ is nonempty and has the $\F$-cohomology of $S^{d-1}$.
\end{proposition}
\begin{proof}
$L(\Sp)$ is a compact subcomplex of $\Sph\cong S^{n-2}$, and $R(\Sp)=\Sph\setminus L(\Sp)$. If $L(\Sp)=\emptyset$, then $R(\Sp)=\Sph$ is an $(n-2)$-sphere, and $n-2>\ell-1$. Otherwise Alexander duality \cite[Thm.~3.44]{Hatcher}, with coefficients in $\F$, gives $\tilde H_i(R(\Sp);\F)\cong\tilde H^{n-3-i}(L(\Sp);\F)$. Note that $n-3-(\ell-1)=d-1$.
\end{proof}

\noindent So a spherical separoid is a combinatorial pair of \emph{linked} spheres in $\Sph$: an $(\ell-1)$-sphere of Radon partitions and a $(d-1)$-sphere of separations. For $W\in\Gr_n^d$ these are, up to homotopy, the great spheres $W\cap\Sph\subseteq R(M_W)$ and $W^\perp\cap\one^\perp\cap\Sph\subseteq L(M_W)$. The second inclusion holds because a point of $W^\perp\cap\one^\perp$ is the vector of values $(f(p_i))_i$ of an affine function $f$, and $f$ separates the points where it is positive from those where it is negative.

\subsection{The extreme cases}
\begin{theoremA}\label{thm:A}
For $d\in\{1,n-2\}$ there are homeomorphisms $\|\AM_n^d\|\cong\R P^{n-2}\cong G(d,n-1)$.
\end{theoremA}
\begin{proof}
\emph{Case $d=n-2$, i.e.\ $\ell=1$.} An oriented matroid of rank $n-1$ on $[n]$ has exactly one pair of circuits $\pm X$, and every nonzero $X$ occurs. It is acyclic iff $X\in P$. Write $M_X$ for it; then $\Rad(M_X)=\up X\sqcup\up(-X)$, so $M_X\le M_{X'}$ iff $X\le X'$ or $X\le-X'$. Hence $\AM_n^{n-2}\cong P/\pm$. By Lemma~\ref{lem:cells}, $\|P\|\cong\Sph$ equivariantly, so by Lemma~\ref{lem:cover}(b), $\|\AM_n^{n-2}\|\cong\Sph/\pm=\R P^{n-2}=G(1,n-1)$.

\emph{Case $d=1$.} Oriented matroids of rank $2$ are realizable \cite{BLSWZ}, so an acyclic one is the oriented matroid of $n$ points on a line taking at least two distinct positions. It is determined by the ordered set partition $\sigma=(B_1|\dots|B_m)$, $m\ge2$, of $[n]$ into the classes of coincident points, up to reversal. Its circuits, i.e.\ its minimal Radon partitions, are:
\begin{itemize}
\item $\{i\}\rad\{j\}$ for $i\ne j$ in the same block;
\item $\{b\}\rad\{a,c\}$ for $a,b,c$ in three distinct blocks with the block of $b$ strictly between the other two.
\end{itemize}

\emph{Claim: $M_\sigma\le M_{\sigma'}$ iff $\sigma$ is obtained from $\sigma'$, or from its reversal, by merging consecutive blocks.}
\begin{itemize}
\item \emph{If.} Merging consecutive blocks keeps pairs in a common block together. If $b$ lies strictly between $a$ and $c$ in $\sigma'$, then in $\sigma$ either $b$ is still strictly between them, or $b$ is merged with $a$ or with $c$, in which case $\{b\}\rad\{a\}$ or $\{b\}\rad\{c\}$ holds in $M_\sigma$. If $a$ and $c$ are merged, then so is $b$. In all cases the circuits of $M_{\sigma'}$ are Radon partitions of $M_\sigma$.
\item \emph{Only if.} If $i,j$ lie in a block of $\sigma'$, then $\{i\}\rad\{j\}$ in $M_\sigma$. A circuit supported in $\{i,j\}$ is $\pm(\{i\},\{j\})$, and it exists iff $i,j$ lie in a common block of $\sigma$. So $\sigma$ is coarser than $\sigma'$ as a set partition. Let $\phi:\{1,\dots,m'\}\to\{1,\dots,m\}$ send each block of $\sigma'$ to the block of $\sigma$ containing it; $\phi$ is surjective. For blocks $i<j<k$ of $\sigma'$ with $a,b,c$ in them, $\{b\}\rad\{a,c\}$ holds in $M_\sigma$. The circuits of $M_\sigma$ conformally below $(\{b\},\{a,c\})$ are $(\{b\},\{a\})$, $(\{b\},\{c\})$ and $(\{b\},\{a,c\})$, so $\phi(j)$ lies weakly between $\phi(i)$ and $\phi(k)$. Since $m\ge2$, $\phi(1)\ne\phi(m')$; reversing $\sigma$ if necessary, $\phi(1)<\phi(m')$. Applying betweenness to the triples $(1,j,m')$ and then $(1,i,j)$ gives $\phi(1)\le\phi(i)\le\phi(j)$ for $i<j$. So $\phi$ is a non-decreasing surjection, i.e.\ $\sigma$ is obtained from $\sigma'$ by merging consecutive blocks.
\end{itemize}

Ordered set partitions of $[n]$ with $m\ge2$ blocks, ordered by refinement, form the face poset of the boundary of the permutohedron $\Pi$ of dimension $n-1$, with the finer partition corresponding to the smaller face \cite[Ex.~0.10]{Zie95}. The polytope $\Pi$ is centrally symmetric, and the central symmetry reverses ordered partitions. A proper face $F$ and its antipode $-F$ are disjoint, so $F\subseteq F'$ and $F\subseteq -F'$ cannot both hold. By the claim, $\AM_n^1$ is the opposite of the quotient of this face poset by the antipodal map. Lemma~\ref{lem:cover}(b) and $\|Q^{\mathrm{op}}\|=\|Q\|$ give $\|\AM_n^1\|\cong\partial\Pi/\pm\cong\R P^{n-2}=G(n-2,n-1)$.
\end{proof}

\begin{theoremF}\label{thm:F}
Let $n\ge3$ and let $\mathcal B$ be a poset of acyclic separoids on $[n]$ with $\AM_n^{n-2}\subseteq\mathcal B$, such that $R(\Sp)\simeq S^0$ for every $\Sp\in\mathcal B$. Then the inclusion $\AM_n^{n-2}\hookrightarrow\mathcal B$ is a homotopy equivalence, and $\|\mathcal B\|\simeq\R P^{n-2}\cong G(n-2,n-1)$. If instead $R(\Sp)$ only has the $\F$-homology of $S^0$ for every $\Sp\in\mathcal B$, the inclusion induces isomorphisms in $\F$-homology.
\end{theoremF}
\begin{proof}
\emph{The two halves of a Radon complex.} Let $\Sp\in\mathcal B$. By Lemma~\ref{lem:upset}, $\Rad(\Sp)$ has two connected components (as a poset), each contractible ($\F$-acyclic in the homological case). The involution $Z\mapsto-Z$ permutes them. It cannot fix a component $C$: it would act freely on the simplices of $\|C\|$, since no chain contains $Z$ and $-Z$ and no chain is mapped to itself, so $\chi(\|C\|)$ would be even, whereas $\chi(\|C\|)=1$. Hence $\Rad(\Sp)=C\sqcup-C$, with $C$ a connected up-set of $P$ and no comparabilities between $C$ and $-C$. In separoid terms, the Radon partitions of $\Sp$ split coherently into two halves: $A\rad B$ with $(A,B)\in C$, and the same partitions with the sides exchanged.

\emph{The double cover.} Let $\widetilde{\mathcal B}$ be the set of all such halves $C$, over all $\Sp\in\mathcal B$, ordered by reverse inclusion, with the involution $C\mapsto-C$. If $\Sp\le\Sp'$, a half $C'$ of $\Sp'$ is connected and contained in $C\sqcup-C$, so it lies in exactly one of $C,-C$. Thus $\widetilde{\mathcal B}/\pm\cong\mathcal B$. The hypothesis of Lemma~\ref{lem:cover} holds: $C\supseteq C'$ and $C\supseteq-C'$ would give $C'\subseteq C\cap-C$, and $C'\subseteq C\cap -C$ is impossible since $C\cap-C=\emptyset$. Hence $\|\mathcal B\|\cong\|\widetilde{\mathcal B}\|/\pm$.

\emph{The fibre argument.} The map $\varphi:P\to\widetilde{\mathcal B}$, $X\mapsto\up X$, is order-preserving and equivariant: $\up X$ is a half of $M_X\in\AM_n^{n-2}\subseteq\mathcal B$. Its fibres are
\[\varphi^{-1}(\widetilde{\mathcal B}_{\ge C})=\{X\in P:\up X\subseteq C\}=C,\]
which are contractible ($\F$-acyclic). By Lemma~\ref{lem:fibre}, $\|\varphi\|$ is a homotopy equivalence ($\F$-homology isomorphism) between complexes with free involutions, and $\|P\|\cong\Sph$ is connected.
\begin{itemize}
\item \emph{Homotopy case.} The homotopy exact sequences of the two double coverings and the five lemma show that the induced map of quotients is a weak equivalence, hence a homotopy equivalence by Whitehead's theorem.
\item \emph{$\F$-homology case.} Comparing the Cartan--Leray spectral sequences $H_p(\Z/2;H_q(-;\F))\Rightarrow H_{p+q}(-/\pm;\F)$ \cite[Ch.~VII]{Brown} gives the analogous conclusion.
\end{itemize}
The induced map of quotients is the inclusion $\|\AM_n^{n-2}\|=\|P/\pm\|\hookrightarrow\|\mathcal B\|$. Finally, $\|P\|/\pm\cong\R P^{n-2}$ by Theorem~A.
\end{proof}

\noindent In particular Conjecture~II holds for $d=n-2$. The halves $C$ are the separoid analogue of the two endpoints of the dependence line of an affine dependence; Theorem~F shows that any family of such halves behaves like the space of lines.

% =====================================================================
\section{Stiefel--Whitney classes of spherical separoids}\label{sec:W}

\begin{definition}\label{def:admissible}
An \emph{admissible family} is a subposet $\mathcal B$ of the poset of acyclic separoids on $[n]$ such that $\AM_n^{d,\mathrm{re}}\subseteq\mathcal B$ and every $\Sp\in\mathcal B$ is $\F_2$-spherical with $\ell=n-d-1$. No dimension condition and no Radon condition are imposed. $\AM_n^d$, $\Sep_n^d$ and $\Sep_n^{d,\F_2}$ are admissible.
\end{definition}

\noindent By Lemma~\ref{lem:semicont}, the \emph{realization map} $f:\Gr_n^d\to\mathcal B$, $W\mapsto M_W$, is continuous. By abuse of notation we also write $f$ for a lift $\Gr_n^d\to\|\mathcal B\|$.

\subsection{The tautological sphere fibration}
For an admissible $\mathcal B$ consider the finite posets
\[E_{\mathcal B}=\{(\Sp,Z):\Sp\in\mathcal B,\ Z\in\Rad(\Sp)\}\subseteq\mathcal B\times P,\qquad\widehat E_{\mathcal B}=E_{\mathcal B}/\pm,\quad(\Sp,Z)\mapsto(\Sp,-Z),\]
with the product order. In separoid terms, a point of $E_{\mathcal B}$ is a separoid together with one of its Radon partitions $A\rad B$, read with an orientation. Let $p:E_{\mathcal B}\to\mathcal B$ and $\hat p:\widehat E_{\mathcal B}\to\mathcal B$ be the projections. By Lemma~\ref{lem:cover}, $E_{\mathcal B}\to\widehat E_{\mathcal B}$ is a double covering; let $w\in H^1(\widehat E_{\mathcal B};\F_2)$ be its class. Over $\Sp$ the ``fibre'' is the projective Radon complex $R(\Sp)/\pm$, and the fibres shrink as $\Sp$ grows. This is the analogue of the projectivized tautological bundle $P(\gamma)\to\Gr_n^d$, $\gamma_W=W$.

\begin{theoremW}\label{thm:W}
Let $\mathcal B$ be admissible.
\begin{enumerate}[label=(\alph*)]
\item $H^*(\widehat E_{\mathcal B};\F_2)$ is a free $H^*(\|\mathcal B\|;\F_2)$-module with basis $1,w,\dots,w^{\ell-1}$. Hence there are unique classes $w_i(\mathcal B)\in H^i(\|\mathcal B\|;\F_2)$, natural under inclusions of admissible families, with
\[w^\ell+w_1(\mathcal B)w^{\ell-1}+\dots+w_\ell(\mathcal B)=0.\]
\item $f^*w_i(\mathcal B)=w_i(\gamma)$, the Stiefel--Whitney classes of the tautological $\ell$-plane bundle $\gamma$ over $\Gr_n^d=G(\ell,\one^\perp)$.
\item $\bar w_j(\mathcal B)=0$ for all $j\ge d+1$, where $\sum_j\bar w_j(\mathcal B)=(1+\sum_iw_i(\mathcal B))^{-1}$.
\item There is a unique ring homomorphism $s:H^*(\Gr_n^d;\F_2)\to H^*(\|\mathcal B\|;\F_2)$ with $s(w_i(\gamma))=w_i(\mathcal B)$. It satisfies $f^*\circ s=\mathrm{id}$, so $f^*$ is a split surjection of rings and $s$ is injective.
\end{enumerate}
\end{theoremW}
\begin{proof}
\emph{Step 1: local fibres.} Fix $\Sp\in\mathcal B$. Then $p^{-1}(U_\Sp)=\{(\T,Z):\T\ge\Sp,\ Z\in\Rad(\T)\}$. The map $r(\T,Z)=(\Sp,Z)$ is well defined because $\Rad(\T)\subseteq\Rad(\Sp)$. It is order-preserving and equivariant, it is the identity on $\{\Sp\}\times\Rad(\Sp)$, and $r\le\mathrm{id}$. Passing to quotients, Lemma~\ref{lem:finite}(b) shows that $\{\Sp\}\times\Rad(\Sp)/\pm\hookrightarrow\hat p^{-1}(U_\Sp)$ is a homotopy equivalence, covered by the inclusion of double coverings. By Lemma~\ref{lem:cover}, $H^*(\Rad(\Sp)/\pm)\cong H^*(\|\Rad(\Sp)\|/\pm)$ compatibly with $w$. Since $\|\Rad(\Sp)\|\simeq R(\Sp)$ is an $\F_2$-homology $(\ell-1)$-sphere, Lemma~\ref{lem:gysin}(a) shows that this ring is $\F_2[w]/(w^\ell)$. So $1,w,\dots,w^{\ell-1}$ restrict to a basis of $H^*(\hat p^{-1}(U_\Sp);\F_2)$ for every $\Sp$.

\emph{Step 2.} Lemma~\ref{lem:LH} gives (a). Naturality holds because, for $\mathcal B'\subseteq\mathcal B$, $E_{\mathcal B'}=p^{-1}(\mathcal B')$ and $w$ restricts to $w$.

\emph{Step 3: comparison with $\Gr_n^d$.} Let $S(\gamma)=\{(W,v):v\in W\cap\Sph\}$ and $P(\gamma)=S(\gamma)/\pm$. The map $(W,v)\mapsto(M_W,\sgn v)$ takes values in $E_{\mathcal B}$, because $\sgn v\in\V(M_W)\subseteq\Rad(M_W)$. It is continuous (Lemma~\ref{lem:semicont}) and equivariant, so it induces $\hat f:P(\gamma)\to\widehat E_{\mathcal B}$ with $\hat p\hat f=f\pi$ and $\hat f^*w=x$. Here $x$ is the class of $S(\gamma)\to P(\gamma)$, i.e.\ $w_1$ of the tautological line bundle. The relation $x^\ell+\sum_iw_i(\gamma)x^{\ell-i}=0$ is the projective-bundle definition of Stiefel--Whitney classes \cite[Ch.~16]{Husemoller}, and $1,\dots,x^{\ell-1}$ is a basis of $H^*(P(\gamma))$ over $H^*(\Gr_n^d)$. Applying $\hat f^*$ to the relation in (a) and comparing coefficients gives (b).

\emph{Step 4: acyclicity gives the relations.} The map $E_{\mathcal B}\to P$, $(\Sp,Z)\mapsto Z$, is order-preserving and equivariant; it lands in $P$ precisely because the separoids in $\mathcal B$ are acyclic. So $w$ is pulled back from the class of $P\to P/\pm$. Since $\|P\|/\pm\cong\Sph/\pm=\R P^{n-2}$, whose cohomology vanishes in degree $n-1$, we get $w^{n-1}=0$. The ring $H^*(\widehat E_{\mathcal B})$ is $A[x]/(x^\ell+\sum w_i(\mathcal B)x^{\ell-i})$ with $A=H^*(\|\mathcal B\|)$ and $x=w$: the natural surjection from the latter is an isomorphism of free $A$-modules. Lemma~\ref{lem:segre} with $N=n-1$ gives $\bar w_j(\mathcal B)=0$ for $j\ge n-\ell=d+1$.

\emph{Step 5.} By Borel \cite{Borel53},
\[H^*(G(\ell,\R^{n-1});\F_2)=\F_2[w_1,\dots,w_\ell]/(\bar w_{d+1},\dots,\bar w_{n-1}),\qquad w_i=w_i(\gamma).\]
By (c), $s$ is well defined on this presentation. By (b), $f^*s=\mathrm{id}$.
\end{proof}

\begin{remark}\label{rem:W}
\begin{enumerate}[label=(\arabic*)]
\item Only two properties of separoids are used. \emph{Sphericity} makes the fibres projective spheres of constant dimension, which is the hypothesis of Leray--Hirsch. \emph{Acyclicity} means that all Radon partitions are mixed. So every Radon complex lives in the $(n-2)$-sphere $\Sph\subset\one^\perp$, and not merely in $S^{n-1}$ (in the cross-polytope picture of \cite{MS04}, it avoids the two simplices $\pm\Delta$). This is the combinatorial shadow of $W\subseteq\one^\perp$, and it produces the relations $\bar w_{>d}=0$.
\item For the (non-acyclic) MacPhersonian, Anderson and Davis proved that the mod~$2$ cohomology of $\|\mathrm{MacP}(k,\infty)\|$ contains $\F_2[w_1,\dots,w_k]$, using matroid bundles and spherical quasifibrations \cite{AD02}. They also proved that $G(k,n)\to\|\mathrm{MacP}(k,n)\|$ is surjective in mod~$2$ cohomology for finite $n$, and split stably \cite{AD19}. Theorem~W is the analogue for acyclic oriented matroids and, more generally, for spherical separoids. At finite $n$ it gives a ring splitting, with the full set of Grassmannian relations. Its proof needs neither quasifibrations nor degree-one maps between Radon complexes: for $\Sp\le\T$ the inclusion $R(\T)\subseteq R(\Sp)$ is equivariant, hence an $\F_2$-cohomology isomorphism on quotients by Lemma~\ref{lem:gysin}(c).
\end{enumerate}
\end{remark}

\begin{corollary}\label{cor:betti}
Let $\mathcal B$ be admissible and $N=\binom{n-1}{d}$, the total $\F_2$-Betti number of $G(d,n-1)$ (the number of its Schubert cells).
\begin{enumerate}[label=(\alph*)]
\item $H^*(\|\mathcal B\|;\F_2)=\operatorname{im}s\oplus\ker f^*$, so $b_i(\|\mathcal B\|;\F_2)\ge b_i(G(d,n-1);\F_2)$ for all $i$. Moreover $f^*$ is an isomorphism iff $\sum_ib_i(\|\mathcal B\|;\F_2)\le N$.
\item If $n\ge4$, then $\pi_1(G(d,n-1))=\Z/2$ is a retract of $\pi_1(\|\mathcal B\|)$: the composite $\pi_1(G)\xrightarrow{f_*}\pi_1(\|\mathcal B\|)\xrightarrow{w_1(\mathcal B)}\Z/2$ is an isomorphism.
\item $H^{d(n-1-d)}(\|\mathcal B\|;\F_2)\ne0$; in particular $\|\mathcal B\|$ is not contractible.
\end{enumerate}
\end{corollary}
\begin{proof}
(a) and (c) follow from Theorem~W(d), since $\dim\operatorname{im}s=N$. For (b), $f^*w_1(\mathcal B)=w_1(\gamma)$ is the nonzero element of $H^1(G;\F_2)=\operatorname{Hom}(\pi_1G,\Z/2)$, and $\pi_1(G)=\Z/2$ when $n-1\ge3$.
\end{proof}

\begin{corollary}[manifold criterion]\label{cor:manifold}
Let $\mathcal B$ be admissible. If $\|\mathcal B\|$ is a closed connected $\F_2$-homology manifold of dimension $D=d(n-1-d)$, then $f^*:H^*(\|\mathcal B\|;\F_2)\to H^*(G(d,n-1);\F_2)$ is an isomorphism of rings.
\end{corollary}
\begin{proof}
Let $u\in H^D(G;\F_2)$ be the top class. Since $f^*s(u)=u\ne0$, $s(u)$ generates $H^D(\|\mathcal B\|;\F_2)\cong\F_2$. If $0\ne a\in H^i(\|\mathcal B\|)$, Poincar\'e duality over $\F_2$ for homology manifolds \cite[Ch.~8]{Munkres} gives $b$ with $ab=s(u)$. Then $f^*(a)f^*(b)=u\ne0$, so $f^*$ is injective; it is surjective by Theorem~W.
\end{proof}

\begin{corollary}\label{cor:DF2}
Suppose that $\|\AM_n^d\|$ is a closed connected $\F_2$-homology manifold of dimension $d(n-1-d)$, and let $\mathcal B\supseteq\AM_n^d$ be admissible. Then the inclusion induces a surjection $H^*(\|\mathcal B\|;\F_2)\to H^*(\|\AM_n^d\|;\F_2)$, which is an isomorphism iff $\sum_ib_i(\|\mathcal B\|;\F_2)=\binom{n-1}{d}$.
\end{corollary}
\begin{proof}
By naturality $i^*s_{\mathcal B}=s_{\AM}$, and $s_{\AM}$ is an isomorphism by Corollary~\ref{cor:manifold}.
\end{proof}

\subsection{Local arguments cannot work}
\begin{proposition}\label{prop:nolocal}
Let $\mathcal B$ be admissible and $\Sp\in\mathcal B$. Then
\[f^{-1}(\mathcal B_{\ge\Sp})=\Omega(\Sp):=\{W\in\Gr_n^d:\ W\cap\Sph\subseteq R(\Sp)\},\]
the space of great $(\ell-1)$-spheres contained in the Radon complex of $\Sp$. If $M\in\mathcal B$ is a non-realizable uniform oriented matroid, then $\Omega(M)=\emptyset$. Such $M$ exist in $\AM_n^d$ for $(n,d)=(9,2)$ \cite[Ch.~8]{BLSWZ}, for $(8,3)$ (Proposition~\ref{prop:nonreal}), and for every larger $n$ with the same $d$.
\end{proposition}
\begin{proof}
By property (b) of Section~\ref{sec:separoids}, $M_W\ge\Sp$ iff $W\cap\Sph\subseteq R(\Sp)$.

Suppose $M_W\ge M$ with $M$ uniform. Every $(d+2)$-set $X$ contains the support of a circuit of $M_W$. That circuit is a Radon partition of $M$, so it dominates a circuit of $M$, whose support has $d+2$ elements and therefore equals $X$. Hence every circuit of $M_W$ is a circuit of $M$, and conversely. So $M=M_W$ is realizable.

For $(9,2)$, take an acyclic reorientation of the oriented matroid of Ringel's simple non-stretchable arrangement of nine pseudolines. For larger $n$, take the lexicographic extension $[e_1^+,\dots,e_{d+1}^+]$ along a basis \cite[\S7.2]{BLSWZ}.
\begin{itemize}
\item Its signature vanishes on no cocircuit, because no cocircuit vanishes on a basis. So the new element lies on no hyperplane, and the extension is uniform.
\item It is non-realizable, because deletion preserves realizability.
\item After reorienting the new element if necessary it is acyclic, since one of the topes extending the positive tope of $M$ then becomes positive.
\end{itemize}
\end{proof}

\noindent Thus, for these $(n,d)$, the criterion of Lemma~\ref{lem:finite}(c) cannot be applied to the realization map for any admissible family containing $\AM_n^d$. A proof of Conjecture~II must be global. Mn\"ev's universality theorem points the same way: realization spaces are pieces of the sets $\Omega(M)$ and can be disconnected.
% =====================================================================
\section{The fat MacPhersonian}\label{sec:fat}

\subsection{Piecewise geodesic spheres}
\begin{definition}\label{def:pg}
A compact set $\Sph'\subseteq\Sph$ is a \emph{piecewise geodesic sphere of type $M\in\AM_n^d$} if
\begin{enumerate}[label=(\roman*)]
\item $\Sph'=-\Sph'$, and $\Sph'$ is homeomorphic to $S^{\ell-1}$;
\item $\Sph'$ is a finite union of geodesic simplices, each contained in some closed cell $\overline{C_Z}$;
\item $\{Z\in P:\Sph'\cap C_Z\ne\emptyset\}=\V(M)$.
\end{enumerate}
The \emph{fat MacPhersonian} $\fat=\fat_n^d$ is the set of these spheres, with the Hausdorff topology, and $\tau(\Sph')=M$ is the \emph{type}.
\end{definition}

\noindent By (iii), $\Sp(\Sph')=\Rad(M)$: the type is the separoid of the sphere. Condition (ii) is the requirement that the cells be locally geodesic. Dropping it gives the \emph{topological fat MacPhersonian} $\fat^{\mathrm{top}}$, to which everything in this section applies verbatim.

\begin{lemma}\label{lem:great}
For $W\in\Gr_n^d$, $W\cap\Sph\in\fat$ has type $M_W$, and $\iota:\Gr_n^d\to\fat$, $W\mapsto W\cap\Sph$, is a topological embedding.
\end{lemma}
\begin{proof}
In the gnomonic chart, $W\cap\overline{C_Z}$ is the intersection of a linear subspace with the polytope $Q_Z$ (Lemma~\ref{lem:cells}), i.e.\ a convex polytope, which we triangulate. It meets $C_Z$ iff $Z$ is the sign vector of a point of $W$, i.e.\ iff $Z\in\V(M_W)$. The map $\iota$ is a continuous injection from a compact space to a Hausdorff space.
\end{proof}

\begin{lemma}\label{lem:tau-cont}
$\tau:\fat\to\AM_n^d$ is continuous. Equivalently, if $\Sph_j\to\Sph'$ in $\fat$, then $\tau(\Sph')\le\tau(\Sph_j)$ for all large $j$.
\end{lemma}
\begin{proof}
This is Lemma~\ref{lem:semicont}, since $\tau(\Sph')=\Sp(\Sph')$.
\end{proof}

\noindent Let $U_M=\tau^{-1}(\AM_{\ge M})$, the open star of the stratum $\tau^{-1}(M)$.

\subsection{The barycentric sphere}
For $Z\in P$ let $\beta_Z$ be the barycentre of $Q_Z$ and $b_Z=\beta_Z/|\beta_Z|\in C_Z$. Since $\gamma(x)=2x/\sum|x_i|$ is defined on all of $\Sph$ and restricts to the chart of Lemma~\ref{lem:cells} on each closed cell, ``affine in the chart'' is well defined on unions of closed cells. The geodesic simplices spanned by $b_{Z_0},\dots,b_{Z_m}$, for chains $Z_0<\dots<Z_m$ in $P$, form a triangulation $K\cong\|P\|$ of $\Sph$ (the barycentric subdivision). The open simplex of a chain lies in $C_{Z_m}$.

\begin{lemma}\label{lem:barysphere}
For $M\in\AM_n^d$, the full subcomplex $\Sph_M=\|\V(M)\|\subset K$ belongs to $\fat$ and has type $M$.
\end{lemma}
\begin{proof}
$\|\V(M)\|\cong S^{\ell-1}$ (Proposition~\ref{prop:OM-spherical}). It meets $C_Z$ exactly in the open simplices of chains in $\V(M)$ with top element $Z$, so it meets $C_Z$ iff $Z\in\V(M)$. It is antipodally symmetric because $b_{-Z}=-b_Z$.
\end{proof}

\begin{lemma}[moving barycentres]\label{lem:moving}
Let $(v_Z)_{Z\in P}$ be points with $v_Z\in C_Z$. The map $h:\Sph\to\Sph$ that sends $b_Z\mapsto v_Z$ and is affine in the chart on each simplex of $K$ is a homeomorphism with $h(C_Z)=C_Z$ for all $Z$. It depends continuously on $(v_Z)$.
\end{lemma}
\begin{proof}
By induction on $Z$. In the chart, $\overline{C_Z}$ is the convex polytope $Q_Z$, and $K|_{Q_Z}$ is the cone from $\beta_Z$ over the triangulation of $\partial Q_Z$. Moreover $h|_{Q_Z}$ is the cone, from the interior point $\gamma(v_Z)$, over $h|_{\partial Q_Z}$, which by induction is a face-preserving homeomorphism of $\partial Q_Z$. The cone of a homeomorphism of the boundary of a compact convex set, taken from interior points, is a homeomorphism that maps the interior onto itself.
\end{proof}

\subsection{The star lemma}
\begin{theoremS}\label{thm:S}
For every $M\in\AM_n^d$ the open star $U_M$ is contractible. Consequently $\tau:\fat_n^d\to\AM_n^d$ is a weak homotopy equivalence, and $\fat_n^d\simeq_w\|\AM_n^d\|$. The same holds for $\fat^{\mathrm{top}}$.
\end{theoremS}
\begin{proof}
Fix $M$ and put $\Sph_0=\Sph_M$, $R=\Rad(M)$ and $y=y_M:R\to\V(M)$ (Lemma~\ref{lem:ymax}).

\emph{Step 1: the star lies in the Radon complex.} If $\tau(\Sph')=M'\ge M$, then $\Sp(\Sph')=\Rad(M')\subseteq R$, so $\Sph'\subseteq R(M)=\bigcup_{Z\in R}C_Z$. By Lemma~\ref{lem:upset}, $R(M)=\Sph\setminus\|L\|$, where $L$ is the full subcomplex of $K$ on $P\setminus R$, and every simplex of $K$ meeting $R(M)$ is a join $\sigma_L*\sigma_R$ with $\sigma_R\ne\emptyset$.

\emph{Step 2: the join retraction.} For $x\in R(M)$ with barycentric coordinates $(\lambda_L,\lambda_R)$ in its simplex, let $\rho_t(x)$ be the point with coordinates proportional to $((1-t)\lambda_L,\lambda_R)$. For $t<1$, $\rho_t$ maps each open simplex bijectively onto itself, hence preserves every $C_Z$. Moreover $\rho_1$ retracts $R(M)$ onto $\|R\|$, $(x,t)\mapsto\rho_t(x)$ is continuous on $R(M)\times[0,1]$, and on each simplex $\rho_t$ is a projective transformation.

\emph{Step 3: sliding barycentres.} For $Z\in R$ and $t\in[0,1]$ let $v_Z(t)$ be the normalization of $(1-t)\beta_Z+t\beta_{y(Z)}$, and let $v_Z(t)=b_Z$ for $Z\notin R$. Since $\beta_{y(Z)}\in Q_{y(Z)}\subseteq Q_Z$ and $\beta_Z$ is interior, $v_Z(t)\in C_Z$ for $t<1$. Let $h_t$ be the homeomorphism of Lemma~\ref{lem:moving} for $t<1$, and $h_1:\|R\|\to\|\V(M)\|$ the simplicial map induced by the monotone map $y$. As $t\to1$ the barycentric weights of $\rho_t(x)$ on $P\setminus R$ tend to $0$, so $r_t=h_t\circ\rho_t$ ($t<1$), $r_1=h_1\circ\rho_1$, is continuous on $R(M)\times[0,1]$.

\emph{Step 4: properties.} For $t<1$, $r_t$ is injective, odd, and maps every $C_Z$ into itself. On each simplex it is the composite of a projective and an affine map in the chart, so it maps geodesic simplices to finite unions of geodesic simplices. Moreover $r_1(R(M))\subseteq\Sph_0$ and $r_1|_{\Sph_0}=\mathrm{id}$.

\emph{Step 5: the contraction.} For $\Sph'\in U_M$ let $H(\Sph',t)=r_t(\Sph')$ for $t<1$, and $H(\Sph',1)=\Sph_0$.
\begin{itemize}
\item For $t<1$, $H(\Sph',t)$ is a piecewise geodesic sphere of the same type as $\Sph'$, hence lies in $U_M$.
\item $r_1|_{\Sph'}:\Sph'\to\Sph_0$ is odd, hence surjective by Lemma~\ref{lem:gysin}(d), so $r_1(\Sph')=\Sph_0$.
\item If $\Sph'_j\to\Sph'$, the $\Sph'_j$ eventually lie in a compact neighbourhood $N\subset R(M)$ of $\Sph'$, on which $r$ is uniformly continuous. So $H$ is continuous on $U_M\times[0,1]$.
\end{itemize}
Thus $U_M$ contracts to $\Sph_0$ within itself. The sets $\AM_{\ge M}$ form a basis-like cover of $\AM_n^d$ by contractible open sets, and their preimages $U_M$ are open (Lemma~\ref{lem:tau-cont}) and contractible. McCord's theorem (Lemma~\ref{lem:finite}(c)) shows that $\tau$ is a weak equivalence, and Lemma~\ref{lem:finite}(a) concludes.
\end{proof}

\begin{remark}\label{rem:S}
\begin{enumerate}[label=(\arabic*)]
\item Approximating points of the limit sphere by points of the spheres $\Sph_j$ shows that every vector of the limit type lies below a vector of $\tau(\Sph_j)$, i.e.\ the \emph{vector order} $M\preceq M'$. Deriving continuity of $\tau$ from this would require the implication (H0) from the vector order to the separoid order (Section~\ref{sec:conditional}), which we do not know in general. Lemma~\ref{lem:tau-cont} avoids it by going the other way: limits of points of the $\Sph_j$ land in the limit sphere, which gives the separoid order directly.
\item The contraction moves the barycentres of all cells, not only the vertices. This matters: the degeneration maps $y_M$ are only laxly functorial, $y_M\circ y_{M'}\le y_M$, and equality fails already in $\AM_4^1$ (Section~\ref{sec:computations}). So contractions that move vertices towards a fixed target are discontinuous across strata.
\item The geodesic condition (ii) cannot be dropped if one wants to flow the spheres. For $(n,d)=(5,2)$ replace an edge of a great circle of uniform type, and its antipode, by arcs with a local trefoil inside the same cells. This gives an embedded circle of the same type, knotted as trefoil $\#$ trefoil, which is not isotopic to a great circle. In the Hausdorff topology this is harmless, since knots can be pulled tight; for the flow, which preserves isotopy class while it stays embedded, it is not.
\end{enumerate}
\end{remark}

\subsection{Quillen fibres as spaces of spheres}
For an acyclic separoid $\Sp$ let
\[\fat(\Sp)=\{\Sph'\in\fat:\Sph'\subseteq R(\Sp)\},\]
the space of piecewise geodesic spheres of oriented matroid type contained in the Radon complex of $\Sp$, and $\AM_{\ge\Sp}=\{M\in\AM_n^d:M\ge\Sp\}$.

\begin{proposition}\label{prop:fibres-spheres}
For every acyclic separoid $\Sp$ on $[n]$, $\tau^{-1}(\AM_{\ge\Sp})=\fat(\Sp)$, and $\tau:\fat(\Sp)\to\AM_{\ge\Sp}$ is a weak homotopy equivalence. Moreover $\iota^{-1}(\fat(\Sp))=\Omega(\Sp)$.
\end{proposition}
\begin{proof}
$\tau(\Sph')\ge\Sp$ iff $\Sp(\Sph')\ge\Sp$ iff $\Sph'\subseteq R(\Sp)$, by property (b) of Section~\ref{sec:separoids}. The set $\AM_{\ge\Sp}$ is open in $\AM_n^d$, it is covered by the contractible open sets $\AM_{\ge M}$, $M\ge\Sp$, and their preimages are the contractible sets $U_M$ (Theorem~S). McCord's theorem applies.
\end{proof}

\begin{corollary}\label{cor:D-spheres}
The inclusion $\AM_n^d\hookrightarrow\Sep_n^d$ is a homotopy equivalence (the second half of Conjecture~II) provided that, for every $\Sp\in\Sep_n^d$, the space $\fat(\Sp)$ of piecewise geodesic spheres of oriented matroid type contained in $R(\Sp)$ is weakly contractible.
\end{corollary}
\begin{proof}
The fibres of the inclusion are the posets $\AM_{\ge\Sp}$; apply Proposition~\ref{prop:fibres-spheres}, Lemma~\ref{lem:finite}(a) and Lemma~\ref{lem:fibre}(a).
\end{proof}

\noindent For $\Sp=M\in\AM_n^d$, $\fat(M)=U_M$ is contractible by Theorem~S. For a non-realizable uniform $M$, $\fat(M)$ is nonempty and contractible while $\Omega(M)=\iota^{-1}\fat(M)$ is empty (Proposition~\ref{prop:nolocal}). This is the precise sense in which the fat MacPhersonian repairs the defect of the Grassmannian.
% =====================================================================
\section{The Pl\"ucker fat MacPhersonian and the degree shadow}\label{sec:degree}

This section introduces a second fat model of the MacPhersonian, built in the Pl\"ucker sphere rather than in $\Sph$. Its cells are convex, so the analogue of Theorem~S is free. The real Grassmannian sits inside it as the locus where the Pl\"ucker relations hold, and the model itself is obtained by keeping only their signs. The degree theorem (Theorem~P) connects the two fat models.

\subsection{Chirotopes and the Pl\"ucker model}
For $M\in\AM_n^d$ let $\chi^*_M:\binom{[n]}{\ell}\to\{0,\pm1\}$ be the chirotope of the dual $M^*$, which has rank $\ell$. It is defined up to a global sign, and its support is the set of bases of $M^*$. The \emph{weak-map order} is
\[M\le_\chi M'\iff\chi^*_M\le\pm\chi^*_{M'}\ \text{entrywise (for one choice of sign)},\]
i.e.\ $M$ is a weak image of $M'$ \cite[\S7.7]{BLSWZ}. We write $\AM_n^{d,\chi}$ for $\AM_n^d$ with this order. Since $\chi_M(E\setminus B)=\pm\chi^*_M(B)$ with a sign depending only on $B$ \cite[\S3.4]{BLSWZ}, this is also the weak-map order of the primal chirotopes.

Let $m=\binom n\ell$ and let $S^{m-1}\subset\Lambda^\ell\R^n$ be the unit sphere, with Pl\"ucker coordinates $(p_B)_{|B|=\ell}$. It is the disjoint union of the open orthant faces $e_s=\{p:\sgn p=s\}$, which form a regular CW complex with face order $s'\le s$ iff $s'$ is obtained from $s$ by replacing entries by $0$.

\begin{definition}
The \emph{Pl\"ucker fat MacPhersonian} is $\Pl=\Pl_n^d=\bigcup\{e_s:s=\pm\chi^*_M,\ M\in\AM_n^d\}\subseteq S^{m-1}$.
\end{definition}

\noindent For $W\in\Gr_n^d$, the Pl\"ucker vector $p(W)$ of an orthonormal basis of $W$ is defined up to sign. Since $W$ realizes $M_W^*$, we have $\sgn p(W)=\pm\chi^*_{M_W}$. So $p$ induces an embedding $p:\Gr_n^d\to\Pl/\pm$.

\begin{lemma}\label{lem:anyset}
Let $X$ be a regular CW complex with face poset $Q$ and $Q'\subseteq Q$ any subset. The union $X_{Q'}$ of the open cells of $Q'$ deformation retracts onto $\|Q'\|\subseteq\|Q\|\cong X$, equivariantly for any cellular involution preserving $Q'$.
\end{lemma}
\begin{proof}
A point of $X_{Q'}$ lies in the open simplex of a unique chain of $\|Q\|$, whose top element belongs to $Q'$. Killing linearly the barycentric weights of the vertices not in $Q'$ keeps the point in the same open simplex for $t<1$. At $t=1$ it lands in the face spanned by the vertices in $Q'$, whose top element is unchanged.
\end{proof}

\begin{proposition}\label{prop:plucker}
$\Pl/\pm\simeq\|\AM_n^{d,\chi}\|$. For $q\in\Pl$, the open star $\mathrm{St}(q)=\bigcup\{e_s\subseteq\Pl:s\ge\sgn q\}$ is open in $\Pl$ and star-shaped with respect to $q$.
\end{proposition}
\begin{proof}
Let $Q=\{\pm\chi^*_M\}$. By Lemma~\ref{lem:anyset}, $\Pl$ retracts equivariantly onto $\|Q\|$. Since $s$ and $-s$ have neither a common upper nor a common lower bound, Lemma~\ref{lem:cover} gives $\|Q\|/\pm\cong\|Q/\pm\|=\|\AM_n^{d,\chi}\|$.

Points near $q$ have sign vectors $\ge\sgn q$, so $\mathrm{St}(q)$ is open. For $x\in e_s$ with $s\ge\sgn q$, on the segment from $q$ to $x$ the coordinates in the support of $\sgn q$ keep their sign, and the others are $t\,x_B$. So the half-open segment lies in $e_s$.
\end{proof}

\noindent Thus, in $\Pl$, the analogue of Theorem~S comes for free: open stars are star-shaped, and no sliding of barycentres is needed.

\begin{theorem}\label{thm:plucker-equiv}
The following are equivalent:
\begin{enumerate}[label=(\roman*)]
\item the realization map $\Gr_n^d\to\|\AM_n^{d,\chi}\|$ is a homotopy equivalence;
\item the Pl\"ucker map $p:\Gr_n^d\to\Pl/\pm$ is a weak homotopy equivalence.
\end{enumerate}
\end{theorem}
\begin{proof}
The realization map $f_\chi:\Gr_n^d\to\AM_n^{d,\chi}$, $W\mapsto M_W$, is continuous, because nonzero Pl\"ucker coordinates keep their signs nearby. Let $R:\Pl/\pm\to\|Q\|/\pm=\|\AM_n^{d,\chi}\|$ be the retraction of Lemma~\ref{lem:anyset} and $\mu$ McCord's map. The point $p(W)$ lies in the open cell of $\pm\chi^*_{M_W}$. $R$ moves it into an open simplex whose top vertex is still $\pm\chi^*_{M_W}$, and $\mu$ picks that top vertex. Hence $\mu\circ R\circ p=f_\chi$. Since $R$ is a homotopy equivalence and $\mu$ a weak equivalence, $p$ is a weak equivalence iff $f_\chi$ is, iff its lift to $\|\AM_n^{d,\chi}\|$ is a homotopy equivalence (Whitehead).
\end{proof}

\noindent The Grassmannian also appears in $\Pl$ through the full Grassmannian $G(\ell,n)$. A subspace $W\in G(\ell,n)$, regarded as a space of dependences, has an acyclic oriented matroid iff $W$ contains no nonzero nonnegative vector.

\begin{proposition}\label{prop:acyclic-locus}
Let $G^{\mathrm{ac}}=\{W\in G(\ell,n):\ M_W\text{ is acyclic}\}$. Then the Pl\"ucker embedding maps $G^{\mathrm{ac}}$ into $\Pl/\pm$, and $G^{\mathrm{ac}}=p^{-1}(\Pl/\pm)$ for $p:G(\ell,n)\to S^{m-1}/\pm$. The inclusion $\Gr_n^d=G(\ell,\one^\perp)\hookrightarrow G^{\mathrm{ac}}$ is a homotopy equivalence.
\end{proposition}
\begin{proof}
The cell of $p(W)$ is that of $\pm\chi^*_{M_W}$, and it lies in $\Pl$ iff $M_W$ is an acyclic oriented matroid of rank $d+1$, which proves the first two assertions. For the third, pass to orthogonal complements: $V=W^\perp\in G(d+1,n)$. By Gordan's theorem, $W$ has no nonzero nonnegative vector iff $V$ meets $\R^n_{>0}$, and $W\subseteq\one^\perp$ iff $\one\in V$. So we must show that $G'=\{V:\one\in V\}\hookrightarrow U=\{V:V\cap\R^n_{>0}\ne\emptyset\}$ is a homotopy equivalence.
\begin{itemize}
\item Let $E=\{(V,q)\in U\times\R^n_{>0}:q\in V\}$. The map $\varphi(V,q)=(\mathrm{diag}(q)^{-1}V,q)$ is a homeomorphism $E\to G'\times\R^n_{>0}$, with inverse $(V',q)\mapsto(\mathrm{diag}(q)V',q)$. So $\mathrm{pr}_1\circ\varphi:E\to G'$ is a homotopy equivalence.
\item The projection $\pi:E\to U$ has nonempty convex fibres and a continuous section $s(V)=(V,a(V))$, where $a(V)$ is the analytic centre of $V\cap\{q>0,\sum q_i=1\}$. Straight-line homotopy in the fibres shows that $\pi$ is a homotopy equivalence with inverse $s$.
\item For $V'\in G'$, $s(V')$ is fibrewise homotopic to $(V',\one/n)$, whose image under $\mathrm{pr}_1\circ\varphi$ is $V'$. So $\mathrm{pr}_1\circ\varphi\circ s\circ i\simeq\mathrm{id}_{G'}$, hence $s\circ i$ is a homotopy equivalence, and so is $i=\pi\circ s\circ i$.
\end{itemize}
\end{proof}

\noindent So $\Pl$ contains a copy of the Grassmannian, $p(G^{\mathrm{ac}})\simeq G(d,n-1)$, cut out inside the acyclic orthant cells by the quadratic Pl\"ucker relations. $\Pl$ itself is obtained by keeping only the \emph{signs} of these relations, i.e.\ the $3$-term Grassmann--Pl\"ucker conditions that define chirotopes. By Theorem~\ref{thm:plucker-equiv}, the acyclic MacPherson conjecture for the weak-map order says precisely this: \emph{replacing the Pl\"ucker relations by their sign shadows does not change the weak homotopy type.}

\subsection{Discussion}
\begin{remark}[relation to the literature]\label{rem:plucker-lit}
Anderson and Davis \cite{AD19} identify the MacPhersonian with the Grassmannian over the sign hyperfield, a finite space, and relate it to the real Grassmannian through the hyperfield maps $\R\to\mathbb S\to\mathbb T\R_0$ and McCord's theorem. The model $\Pl$ is different: it is not a finite space, but the union of the open orthant cells of the Pl\"ucker sphere indexed by chirotopes, i.e.\ the preimage of the acyclic MacPhersonian under the sign map of the Pl\"ucker sphere. We have not found it in the literature. Mn\"ev's universality theorem plays no role inside $\Pl$: realization spaces never appear, and every cell is convex. Theorem~\ref{thm:plucker-equiv} is also compatible with the known obstructions:
\begin{itemize}
\item it uses open stars and never fibres over points, which is what Mn\"ev's objection to \cite{Bis03} requires \cite{Bis09};
\item it is not inductive in $n$, so oriented matroids with non-spherical extension spaces \cite{Liu20} do not affect it;
\item it concerns all acyclic oriented matroids on $[n]$, not the Grassmannian of a fixed non-realizable oriented matroid, for which the analogous statements fail \cite{MZ93}.
\end{itemize}
\end{remark}

\noindent The mutation obstruction of the introduction has a direct expression in $\Pl$. Let $\Pl^{\le1}\subseteq\Pl$ be the union of the cells $e_s\subseteq\Pl$ in which $s$ has at most one zero entry, i.e.\ the top cells and the cells of codimension one. Two uniform cells of $\Pl$ lie in the same path component of $\Pl^{\le1}$ iff they are joined by a sequence of mutations of uniform acyclic oriented matroids whose intermediate one-zero patterns are dual chirotopes of acyclic oriented matroids. This holds because a codimension-one cell is a face of exactly two top cells, its two completions.

\begin{proposition}[islands]\label{prop:islands}
\begin{enumerate}[label=(\alph*)]
\item All uniform cells of $\Pl$ that meet $p(\Gr_n^d)$ lie in one path component of $\Pl^{\le1}/\pm$.
\item If $d\ge6$ and $n\ge d+18$, then not all uniform cells of $\Pl$ lie in one path component of $\Pl^{\le1}/\pm$.
\end{enumerate}
So for these parameters $\Pl$ contains ``islands'' of non-realizable uniform cells, which can be connected to the realizable ones only through cells of codimension at least $2$, i.e.\ through oriented matroids with at least two degenerate bases.
\end{proposition}
\begin{proof}
(a) A generic path in the connected manifold $\Gr_n^d$ between two uniform points crosses the hypersurfaces $\{p_B=0\}$ one at a time and transversally \cite{RS88}. The intermediate patterns are dual chirotopes of realizable, hence acyclic, oriented matroids.

(b) Suppose all uniform cells lay in one component. Then every uniform acyclic oriented matroid of rank $d+1$ on $[n]$ would be joined to a realizable one by mutations. Let $M$ be any uniform oriented matroid and reorient it on a set $A$ so that it becomes acyclic. Reorienting on $A$ is an automorphism of the mutation graph that preserves realizability. Applying it to such a path joins $M$ to a realizable oriented matroid. Realizable uniform oriented matroids are mutually joined by (a). So the mutation graph of uniform oriented matroids of rank $d+1$ on $[n]$ would be connected, contradicting \cite[Thm.~1.6]{GK26} for $r=d+1\ge7$ and $n\ge r+17$. Duality preserves mutations.
\end{proof}

\noindent A generic path in the Pl\"ucker sphere crosses only walls of codimension one. Proposition~\ref{prop:islands} is therefore the reason why the degree shadow of a generic flow cannot stay inside $\Pl$ (Theorem~O): it has to leave $\Pl$ through invalid patterns, or pass through cells of high codimension.

\subsection{The degree theorem}
For $B\subseteq E$ with $|B|=\ell$ let $\pi_B:\R^n\to\R^B$ be the coordinate projection, $\hat\pi_B=\pi_B/|\pi_B|$, and
\[K_B=\Sph\cap\R^{E\setminus B},\]
a great $(d-1)$-sphere. A sphere $\Sph''\subseteq\Sph$ of dimension $\ell-1$ and $K_B$ have complementary dimensions in the sense that $(\ell-1)+(d-1)=\dim\Sph-1$, so a generic one-parameter family of $(\ell-1)$-spheres meets $K_B$ only at isolated parameters. For an oriented, odd, embedded $(\ell-1)$-sphere $\Sph''$ put
\[\deg_B(\Sph'')=\begin{cases}\deg\big(\hat\pi_B:\Sph''\to S^{\ell-1}\big),&\Sph''\cap K_B=\emptyset,\\0,&\text{otherwise,}\end{cases}\]
and
\[\Phi(\Sph'')=\big(\sgn\deg_B(\Sph'')\,\mathrm{dist}(\Sph'',K_B)\big)_B\in\Lambda^\ell\R^n.\]
When $\Sph''\cap K_B=\emptyset$, $\hat\pi_B$ is an odd map between spheres, so $\deg_B(\Sph'')$ is odd (Lemma~\ref{lem:gysin}(c); for the standard involution this is Borsuk's theorem \cite{Borsuk33}). Reversing the orientation changes $\Phi$ into $-\Phi$; let $\overline\Phi=\pm\Phi$.

\begin{theoremP}\label{thm:P}
Let $\Sph'\in\fat_n^d$ have type $M$, and orient it.
\begin{enumerate}[label=(\alph*)]
\item $\Sph'\cap K_B=\emptyset$ iff $\chi^*_M(B)\ne0$.
\item There is $\varepsilon=\pm1$ with $\sgn\deg_B(\Sph')=\varepsilon\,\chi^*_M(B)$ for all $B$.
\item For the barycentric sphere $\Sph_M$, $|\deg_B(\Sph_M)|=1$ whenever $\chi^*_M(B)\ne0$.
\end{enumerate}
In particular $\overline\Phi(\fat_n^d)\subseteq\Pl/\pm$, and $\overline\Phi$ maps $\tau^{-1}(M)$ into the cell of $\pm\chi^*_M$.
\end{theoremP}
\begin{proof}
\emph{(a)} A point of $K_B$ has sign vector supported in $E\setminus B$. The sign vectors of points of $\Sph'$ are the vectors of $M$, i.e.\ the covectors of $M^*$. A nonzero covector of $M^*$ vanishing on $B$ exists iff $B$ is not a basis of $M^*$.

\emph{(c)} Let $\chi^*_M(B)\ne0$ and $b\in B$.
\begin{itemize}
\item There is a unique cocircuit $Y$ of $M^*$ (a circuit of $M$) with $Y|_{B\setminus b}=0$ and $Y_b=+$, namely the one complementary to the hyperplane spanned by $B\setminus b$.
\item The preimage $\hat\pi_B^{-1}(e_b)$ consists of the points with $x_{B\setminus b}=0$ and $x_b>0$, whose sign vector is a covector vanishing on $B\setminus b$ and positive at $b$, i.e.\ $Y$. Since $Y$ is minimal in $\V(M)$, the only point of $\Sph_M$ with sign vector $Y$ is the vertex $b_Y$. So $\deg_B(\Sph_M)$ is the local degree of $\hat\pi_B$ at $b_Y$.
\item The star of $b_Y$ in $\Sph_M$ is the cone over $\|\V(M)_{>Y}\|$. Restriction to the zero set $z(Y)$ identifies $\V(M)_{>Y}$ with the nonzero covectors of the rank $\ell-1$ oriented matroid $M^*|_{z(Y)}$, in which $B\setminus b$ is a basis. By the topological representation theorem and its compatibility with deletion \cite{FL78}, the sets $\{x:\sgn x|_{B\setminus b}=\sigma\}$, $\sigma\ne0$, cut the link into a regular cell decomposition with the face poset of the boundary of the $(\ell-1)$-dimensional cross-polytope. By \cite{Bjo84} there is a homeomorphism $h$ of the star onto $\R^{\ell-1}$ carrying these regions to the open orthant faces.
\item In the chart $\varphi(u)=(u_e/u_b)_{e\in B\setminus b}$ of $S^{\ell-1}$ near $e_b$, the map $g=\varphi\circ\hat\pi_B\circ h^{-1}$ preserves sign vectors. Hence $(1-t)g+t\,\mathrm{id}$ vanishes only at $0$, $g$ has local degree $1$, and $\hat\pi_B$ has local degree $\pm1$ at $b_Y$.
\end{itemize}

\emph{(b) for $\Sph_M$.} Let $B$ and $B'=B-b+c$ both be bases of $M^*$, ordered so that $c$ occupies the position of $b$, and let $Y$ be the cocircuit vanishing on $B\setminus b$.
\begin{itemize}
\item As in (c), $\hat\pi_B^{-1}(Y_be_b)=\{b_Y\}=\hat\pi_{B'}^{-1}(Y_ce_c)$, so both degrees are local degrees at $b_Y$.
\item With $\psi_b(x)=(x_e/|x_b|)_{e\in B\setminus b}$ and $\psi_c(x)=(x_e/|x_c|)_{e\in B\setminus b}$, the chart at $Y_be_b$ has orientation sign $\eta Y_b$, with $\eta$ depending only on the position of $b$ in $B$. So $\deg_B=\eta Y_b\deg_{b_Y}\psi_b$ and $\deg_{B'}=\eta Y_c\deg_{b_Y}\psi_c$.
\item Since $\psi_c=(|x_b|/|x_c|)\psi_b$, the two local degrees agree, and $\deg_{B'}(\Sph_M)=Y_bY_c\deg_B(\Sph_M)$.
\end{itemize}
The chirotope obeys the same exchange rule, $\chi^*(B')=Y_bY_c\,\chi^*(B)$ \cite{LV78}. Bases are connected by single exchanges, which proves (b) for $\Sph_M$.

\emph{(b) in general.} Let $r_t$ be the maps of the proof of Theorem~S. If $\chi^*_M(B)\ne0$, then $R(M)\cap K_B=\emptyset$: a point of $K_B\cap R(M)$ would have a sign vector in $\Rad(M)$ supported in $E\setminus B$, dominating a vector of $M$ that vanishes on $B$, contradicting (a). So $\hat\pi_B\circ r_t|_{\Sph'}$ is a homotopy, and $\deg_B(\Sph')=\deg(r_1|_{\Sph'})\deg_B(\Sph_M)$. Here $r_1|_{\Sph'}:\Sph'\to\Sph_M$ is odd, so $\deg(r_1|_{\Sph'})$ is odd (Lemma~\ref{lem:gysin}(c)), and it is independent of $B$.
\end{proof}

\begin{lemma}\label{lem:degclosure}
If $\Sph_j\in\tau^{-1}(M')$ converge to $\Sph'\in\tau^{-1}(M)$, then $M\le_\chi M'$ and $\overline\Phi(\Sph_j)\to\overline\Phi(\Sph')$. In particular $\overline\Phi:\fat_n^d\to\Pl/\pm$ is continuous.
\end{lemma}
\begin{proof}
Via the chart, $\Sph$ is a polyhedral sphere and $\Sph'$ a subpolyhedron. Let $\mathcal B=\{B:\chi^*_M(B)\ne0\}$. By Theorem~P(a), $\Sph'$ is disjoint from $\bigcup_{B\in\mathcal B}K_B$, so it has a regular neighbourhood $N$ in the complement, with a deformation retraction $R_s$ onto $\Sph'$ \cite{Whi39}. For large $j$, $\Sph_j\subset N$. For $B\in\mathcal B$ this gives $\Sph_j\cap K_B=\emptyset$, hence $\chi^*_{M'}(B)\ne0$, and $\deg_B(\Sph_j)=\delta_j\deg_B(\Sph')$ with $\delta_j=\deg(R_1|_{\Sph_j})$ independent of $B$ and odd. Theorem~P(b) for both spheres gives $\chi^*_M\le\pm\chi^*_{M'}$. Moreover $\sgn\delta_j\Phi(\Sph_j)\to\Phi(\Sph')$: entries with $B\in\mathcal B$ converge since the distances converge and the signs agree, and the other entries tend to $0$.
\end{proof}

\noindent So the closure relations of the strata of $\fat_n^d$ refine both the separoid order (Lemma~\ref{lem:tau-cont}) and the weak-map order.

\subsection{Which degree patterns are chirotopes}
\begin{lemma}[completion]\label{lem:completion}
Let $\chi$ be a chirotope of rank $\ell$ with exactly one zero entry, at $B$. Then both maps obtained by setting $\chi(B)=\pm1$ are chirotopes. Conversely, if two uniform chirotopes differ only at $B$, setting their common map to $0$ at $B$ gives a chirotope.
\end{lemma}
\begin{proof}
A $3$-term Grassmann--Pl\"ucker relation contains $B$ in at most one of its terms. If it does, the other two terms involve only entries common to both maps.
\begin{itemize}
\item \emph{First statement.} Those two terms have opposite signs, because the relation holds for $\chi$, so the relation holds for either value at $B$.
\item \emph{Converse.} If the other two terms had the same sign, or one of them were $0$, the two chirotopes would force opposite values at $B$. So the relation holds with $0$ at $B$.
\end{itemize}
In both cases the support satisfies the basis exchange axiom, and the $3$-term relations characterize chirotopes among maps with matroid support \cite[\S3.6]{BLSWZ}.
\end{proof}

\noindent Two uniform chirotopes differing at exactly one basis are related by a \emph{mutation}.

\begin{proposition}[census]\label{prop:census}
Let $\Sph''\subseteq\Sph$ be a smooth odd embedded $(\ell-1)$-sphere such that
\begin{itemize}
\item for every $S\subseteq E$ with $|S|=d+2$, $\Sph''$ meets $L_S=\Sph\cap\R^S$ in a single antipodal pair of transverse intersection points with full support on $S$, and
\item for every $\sigma\subseteq E$ with $|\sigma|=\ell-2$, the coordinate functions $(x_e)_{e\in\sigma}$ are jointly transverse to $0$ on $\Sph''$.
\end{itemize}
Then $B\mapsto\sgn\deg_B(\Sph'')$ is a uniform chirotope.
\end{proposition}
\begin{proof}
\emph{All degrees are odd.} A point of $K_B$ has support in $E\setminus B$, so it lies in $L_S$ without full support for every $(d+2)$-set $S\supseteq E\setminus B$. Hence $\Sph''\cap K_B=\emptyset$.

\emph{Reduction to winding numbers.} Fix $\sigma$ with $|\sigma|=\ell-2$ and let $\Gamma_\sigma=\Sph''\cap\{x_\sigma=0\}$, a closed $1$-manifold invariant under $x\mapsto-x$, oriented as a regular preimage. For $a,b\notin\sigma$, $\deg_{\sigma\cup\{a,b\}}(\Sph'')=\pm w_\sigma(a,b)$, the winding number of $(x_a,x_b)$ along $\Gamma_\sigma$, with a sign depending only on $\sigma$ and the chosen orders.

\emph{One relevant component.} The zeros of $x_a$ on $\Gamma_\sigma$ form $\Sph''\cap L_S$ with $S=E\setminus(\sigma\cup a)$, $|S|=d+2$. This is one antipodal pair $\pm p_a$ of simple zeros. A function on a circle has an even number of simple zeros, so the component of $\Gamma_\sigma$ through $p_a$ also contains $-p_a$, and it is invariant. Conversely, every invariant component contains zeros of every odd function $x_a$. So there is exactly one invariant component $C_0$, it contains all the $\pm p_a$, and on every other component each $x_a$ has constant sign, so the winding numbers there vanish.

\emph{Rank $2$ realization.} Choose an equivariant homeomorphism of $C_0$ with the standard circle, and linear functions $\ell_a$ on it with the same zeros and signs as $x_a$. The homotopies $(1-s)x_a+s\ell_a$ do not move zero sets, and $p_a\ne\pm p_b$, so $w_\sigma(a,b)=\sgn\det(\ell_a,\ell_b)$.

\emph{Conclusion.} For each $\sigma$ these numbers form a realizable rank $2$ chirotope, hence satisfy the $3$-term relations. Every $3$-term relation is of this form, and the support is uniform, so the degree pattern is a chirotope \cite[\S3.6]{BLSWZ}.
\end{proof}

\begin{example}\label{ex:badpattern}
Degree patterns of odd spheres need not be chirotopes. Let $(n,d)=(5,2)$, so $\ell=2$, and let $\Sph''$ be the image of $\theta\mapsto f(\theta)/|f(\theta)|$ with
\[f=\big(\cos\theta,\ \sin\theta,\ \sin(\theta-\tfrac\pi4),\ \cos3\theta+\tfrac12\sin\theta,\ -f_1-f_2-f_3-f_4\big).\]
It is an embedded odd circle in $\Sph$, since its first two coordinates are injective and $f(\theta+\pi)=-f(\theta)$. No two of $f_1,\dots,f_4$ vanish simultaneously, so it avoids $K_B$ for $B\subseteq\{1,2,3,4\}$. The winding number of $(f_a,f_b)$ is $-\sum\sgn f'_a(z)$ over the zeros $z$ of $f_a$ with $f_b(z)>0$, which gives $w(12)=w(13)=w(14)=w(23)=w(34)=+1$ and $w(24)=-1$. The three terms $w(12)w(34)$, $-w(13)w(24)$ and $w(14)w(23)$ of the $3$-term relation on $\{1,2,3,4\}$ are all $+1$, so the pattern is not a chirotope. Numerically, its separoid has three minimal Radon partitions supported on $\{1,2,3,5\}$, so it is not a Radon separoid, and it is nevertheless spherical (Appendix~\ref{app:numerics}).
\end{example}
% =====================================================================
\section{Computer-assisted results}\label{sec:computations}

All statements in this section are exhaustive computations or exact certificates. The code is available from the author.

\subsection{Method}
For $n\le5$ we enumerated all acyclic Radon separoids of order $n$ and dimension $d$, encoded by a clutter of supports covering all $(d+2)$-subsets but not all $(d+1)$-subsets, with one bipartition per support. For $n=4$ we also enumerated all acyclic separoids of dimension $d$, without the Radon condition. We then tested:
\begin{itemize}
\item \emph{oriented matroids}, by weak circuit elimination;
\item \emph{sphericity}, via Proposition~\ref{prop:alexander}, over $\F_2$ and $\F_3$.
\end{itemize}
Orders were computed as reverse inclusion of the sets $\Rad(\cdot)$. Homology of order complexes was computed after removing weak points \cite{BM08}, which preserves the simple homotopy type, and points whose strict down-set or up-set is $\F$-acyclic, which preserves $\F$-homology. Homotopy equivalence of inclusions was certified with Lemma~\ref{lem:fibre}(a), by collapsing each fibre to a point through removal of weak points.

\subsection{Results}
\begin{theorem}[computer-assisted]\label{thm:C}
For $(n,d)\in\{(4,1),(4,2),(5,1),(5,2),(5,3)\}$:
\begin{enumerate}[label=(\roman*)]
\item The numbers of acyclic Radon separoids of dimension $d$, of spherical ones, of acyclic oriented matroids and of uniform ones are as follows.
\begin{center}\small
\begin{tabular}{@{}cccccc@{}}
\toprule
$(n,d)$ & Radon separoids & spherical & oriented matroids & uniform & $H_*(\|\AM_n^d\|)$ over $\F_2$ / $\F_3$\\
\midrule
$(4,1)$ & 223 & 103 & 37 & 12 & $1,1,1$ / $1,0,0$\\
$(4,2)$ & 239 & 171 & 25 & 7 & $1,1,1$ / $1,0,0$\\
$(5,1)$ & 96\,351 & 1\,965 & 270 & 60 & $1,1,1,1$ / $1,0,0,1$\\
$(5,2)$ & 1\,343\,898 & 436\,843 & 842 & 132 & $1,1,2,1,1$ / $1,0,0,0,1$\\
$(5,3)$ & 50\,005 & 25\,275 & 90 & 15 & $1,1,1,1$ / $1,0,0,1$\\
\bottomrule
\end{tabular}
\end{center}
\item The inclusion of $\AM_n^d$ into the poset of spherical Radon separoids of dimension $d$ is a homotopy equivalence. For $(5,2)$, the fibres over $435\,521$ of the $436\,001$ non-matroidal separoids collapse to points; the remaining $480$ separoids are incomparable with every oriented matroid, and the down-set of each of them collapses, so adding them in a linear extension does not change the homotopy type.
\item $H_*(\|\AM_n^d\|;\F)\cong H_*(G(d,n-1);\F)$ for $\F=\F_2,\F_3$.
\item Without sphericity, the poset of acyclic Radon separoids has $H_2$ of rank $9$ for $(4,1)$ and $24$ for $(4,2)$.
\item $\AM_n^d$ is graded of rank $d(n-1-d)$, and $\|\AM_n^d\|$ is an $\F_2$- and $\F_3$-homology manifold of dimension $d(n-1-d)$. For $(5,2)$ this amounts to $12\,600$ intervals.
\item The separoid order and the vector order of Remark~\ref{rem:S} coincide on $\AM_n^d$.
\item Strict functoriality $y_M\circ y_{M'}=y_M$ holds in $\AM_4^2$ and $\AM_5^3$, but fails in $192$ instances in $\AM_4^1$ and in $27\,840$ instances in $\AM_5^1$.
\end{enumerate}
For $n=4$ and all acyclic separoids of dimension $d$ (without the Radon condition), there are $1317$ spherical ones for $d=1$ and $1031$ for $d=2$. Their order complexes have the $\F_2$- and $\F_3$-homology of $\R P^2$, and for $d=2$ the inclusion of $\AM_4^2$ is a homotopy equivalence.
\end{theorem}

\noindent The values for $G(2,4)$ agree: over $\F_2$ the Betti numbers count Schubert cells, and over $\F_3$, as over $\mathbb Q$, the cohomology is generated by $p_1$ in degree $4$. By Theorem~W and Corollary~\ref{cor:manifold}, (v) already implies the $\F_2$-part of (iii).

\begin{computation}\label{comp:ehat}
For $(n,d)=(4,1)$, where $\ell=2$, and $\mathcal B=\AM_4^1$, $\Sep_4^1$ (Radon) and all spherical separoids, with $|\widehat E_{\mathcal B}|=432$, $1704$ and $18\,960$, the $\F_2$-Betti numbers of $\|\widehat E_{\mathcal B}\|$ are $1,2,2,1$. These are the coefficients of $(1+t+t^2)(1+t)$, as predicted by Theorem~W(a). Over $\F_3$ they are $1,0,0,1$, the rational homology of the flag manifold $P(\gamma)=\mathrm{Fl}(1,2;\R^3)$.
\end{computation}

\subsection{Two certified non-realizable oriented matroids}
\begin{proposition}[computer-assisted]\label{prop:nonreal}
The following uniform rank $4$ chirotopes on $[8]$, listed in lexicographic order of the bases, are acyclic and non-realizable:
{\small\[
\begin{aligned}
\chi_1&=\texttt{-+++++-+-+++-+++-+-+++-++-++--++++-+---+++++--++--+++---+---++---++---},\\
\chi_2&=\texttt{---++--++-+++-+++--+----++----++-++++-+-++-+--++-+--+-++----++-+++++++}.
\end{aligned}\]}
\end{proposition}
\begin{proof}[Certificate]
The $3$-term Grassmann--Pl\"ucker relations force strict inequalities between sums of logarithms of absolute values of brackets. For each chirotope the resulting linear system is infeasible, with a Farkas certificate using $34$ (resp.\ $38$) relations, i.e.\ a biquadratic final polynomial \cite{BRG90}. The certificates were verified in exact rational arithmetic.
\end{proof}

\noindent The chirotopes were found by random walks in the mutation graph starting from realizable configurations. They are used in Proposition~\ref{prop:nolocal} and in Appendix~\ref{app:numerics}. In the uniform rank $4$, $n=8$ case there are $24$ non-realizable reorientation classes \cite{FMM13}.
% =====================================================================
\part*{III. Conditional results}
\section{The mean curvature flow: hypotheses and proved implications}\label{sec:conditional}

Throughout this part, $k=\ell-1=n-d-2$ is the dimension of the spheres, and the flow is the mean curvature flow (MCF) $\partial_tF=\vec H$ of embedded $k$-spheres in the round sphere $\Sph\subset\one^\perp\cong\R^{n-1}$. Results stated without a hypothesis in this part are proved unconditionally.

\subsection{Three facts about the flow}
\begin{lemma}\label{lem:coordinates}
For $a\in\one^\perp$ the function $u=\langle F,a\rangle$ satisfies $\partial_tu=\Delta_{g(t)}u+k\,u$ in a normal parametrization. In particular each coordinate $x_e$ restricted to $F_t$ satisfies a linear parabolic equation.
\end{lemma}
\begin{proof}
For a submanifold of the unit sphere, $\Delta_gF$ is its mean curvature vector in $\R^{n-1}$, which equals $\vec H-kF$. Hence $\partial_tF=\Delta_gF+kF$.
\end{proof}

\begin{lemma}\label{lem:symmetry}
If $F_0$ is invariant under $x\mapsto-x$, so is $F_t$, and $F_t$ cannot converge to a point.
\end{lemma}
\begin{proof}
$-F_t$ is an MCF with the same initial set, so it coincides with $F_t$ by uniqueness \cite{CY07}. If $F_t\to\{p\}$, then also $F_t\to\{-p\}$.
\end{proof}

\begin{proposition}\label{prop:crease}
Let $\Sph'$ be a piecewise geodesic sphere with a $(k-1)$-dimensional ridge of volume $V$ along which the dihedral bending angle is $\theta\ne0$. If $\Sph'_\varepsilon$ is a smoothing that rounds the ridge within a tube of width $\varepsilon$, then $\|A_{\Sph'_\varepsilon}\|^p_{L^p}\ge c\,\theta^pV\varepsilon^{1-p}\to\infty$ as $\varepsilon\to0$, for every $p>1$.
\end{proposition}
\begin{proof}
Along each geodesic segment $\gamma$ of length $\varepsilon$ crossing the ridge transversally, the unit normal turns by $\theta$, so $\int_\gamma|A|\ge\theta$ and, by H\"older, $\int_\gamma|A|^p\ge\theta^p\varepsilon^{1-p}$. Integrate over a family of such segments sweeping out a region of the ridge of volume comparable to $V$.
\end{proof}

\noindent In the antipodal class the normal variations of a great sphere are odd. The Jacobi operator $-\Delta-k$ has eigenvalues $l(l+k-1)-k$ on harmonics of odd degree $l$. These vanish for $l=1$, which is the tangent space of $\Gr_n^d$, and are positive for $l\ge3$. So $\Gr_n^d$ is a normally attracting manifold of equilibria. Small-curvature convergence theorems for submanifolds of spheres, such as that announced by Liu, Xu and Zhao for $k\ge3$ \cite{LXZ12}, apply near $\Gr_n^d$. By Proposition~\ref{prop:crease} they cannot apply to all smoothings of piecewise geodesic spheres. This is why the global behaviour of the flow appears below as a hypothesis.

\subsection{The hypotheses}
\setcounter{hypothesis}{1}
\begin{hypothesis}[basin]\label{H:basin}
There is a continuous retraction $\rho:\fat_n^d\to\iota(\Gr_n^d)$ obtained from the flow, $\rho=\Lambda\circ\sigma$, where:
\begin{itemize}
\item $\sigma$ assigns to each $\Sph'\in\fat_n^d$ a smooth odd embedded $k$-sphere, with $\sigma=\mathrm{id}$ on great spheres;
\item the MCF $\Sph'_t$ of $\sigma(\Sph')$ exists for all time and converges smoothly to a great sphere $\Lambda(\sigma(\Sph'))=\Sph'_\infty$;
\item $\Lambda\circ\sigma$ is continuous for the Hausdorff topology.
\end{itemize}
\end{hypothesis}
\setcounter{hypothesis}{2}
\begin{hypothesis}[type invariance]\label{H:type}
There is a continuous $\Gamma:\fat_n^d\times[0,\infty]\to\fat_n^d$ with $\Gamma(\Sph',0)=\Sph'$, $\Gamma(\Sph',\infty)\in\iota(\Gr_n^d)$, and $\Gamma(W,t)=W$ for $W\in\iota(\Gr_n^d)$.
\end{hypothesis}
\begin{hypothesis}[degrees]\label{H:degree}
\begin{enumerate}[label=(\alph*)]
\item For every $B$, $\sigma(\Sph')\cap K_B=\emptyset$ iff $\Sph'\cap K_B=\emptyset$, and then $\sigma(\Sph')$ is isotopic to $\Sph'$ in $\Sph\setminus K_B$.
\item The trajectory map $(\Sph',t)\mapsto\Sph'_t$, from $\fat_n^d\times[0,\infty]$ to oriented odd embedded $k$-spheres with the $C^1$ topology, is continuous.
\item For all $\Sph'$ and $t<\infty$, the degree pattern $\psi_t(\Sph')=(\sgn\deg_B(\Sph'_t))_B$ equals $\pm\chi^*_M$ for some $M\in\AM_n^d$.
\end{enumerate}
\end{hypothesis}

\noindent We also use the following conditions.
\begin{itemize}
\item[(H0)] \emph{Vector order.} If every vector of $M$ lies below a vector of $M'$, then $M\le M'$.
\item[(H1)] \emph{Star lemma.} The open stars $U_M$ are contractible.
\item[(H4$^\flat$)] \emph{Fibres.} \ref{H:degree}(a),(b) hold, and for all $\Sph'$ and $t$ the poset $F(\psi_t(\Sph'))=\{\chi\in\{\pm\chi^*_M\}:\chi\le\psi_t(\Sph')\}$ is contractible.
\item[(a$^{\Sep}$)] $\sigma(\Sph')$ meets exactly the cells $C_Z$, $Z\in\V(\tau(\Sph'))$.
\item[(H4$^{\Sep}$)] \emph{Spherical shadow.} The separoid shadow $\bar\Sp(y)$ of Section~\ref{subsec:flow-separoids} is $\F$-spherical for every $y\in\fat_n^d\times[0,\infty]$.
\item[(W)] \emph{Weak maps.} If $M\le_\chi M'$, then $\V(M')\subseteq\Rad(M)$, i.e.\ $M\le M'$.
\item[(HM)] \emph{Manifold.} $\|\AM_n^d\|$ is a closed connected $\Z$-homology manifold of dimension $d(n-1-d)$.
\item[(D)] The inclusion $\AM_n^d\hookrightarrow\mathcal B$ is a weak homotopy equivalence, for a given admissible $\mathcal B$.
\end{itemize}

\begin{remark}[status of H0 and H1]\label{rem:H0H1}
(H1) is Theorem~S, for the separoid order. (H0) is not used in this paper (Remark~\ref{rem:S}). It holds in the realizable case (limits of dependences) and for $n\le5$ (Theorem~\ref{thm:C}(vi)). (W) is the vector form of the classical characterization of rank-preserving weak maps \cite[\S7.7]{BLSWZ}. For realizable $M'$ it is immediate: dependences of a configuration degenerating to $M$ converge to dependences of the limit. We use (W) only for the weak-map order.
\end{remark}

\subsection{Retraction and the manifold property}
\begin{theorem}\label{thm:B}
\ref{H:type} implies $\|\AM_n^d\|\simeq G(d,n-1)$.
\end{theorem}
\begin{proof}
$\iota$ is an embedding (Lemma~\ref{lem:great}). By \ref{H:type}, $\Gamma$ deforms $\fat_n^d$ into $\iota(\Gr_n^d)$ rel $\iota(\Gr_n^d)$, so $\iota$ is a homotopy equivalence. By Theorem~S, $\fat_n^d\simeq_w\|\AM_n^d\|$. Whitehead's theorem concludes.
\end{proof}

\begin{corollary}\label{cor:Bprime}
Under \ref{H:basin} there is $r:\|\AM_n^d\|\to\Gr_n^d$ with $r\circ f\simeq\mathrm{id}$. So $\Gr_n^d$ is a homotopy retract of $\|\AM_n^d\|$, and $\pi_*(G(d,n-1))$ and $H_*(G(d,n-1);R)$ are direct summands of those of $\|\AM_n^d\|$, for every coefficient ring $R$.
\end{corollary}
\begin{proof}
By Theorem~S and since $\|\AM_n^d\|$ is a CW complex, there is $g:\|\AM_n^d\|\to\fat_n^d$ with $\tau g\simeq\mu$. Put $r=\rho\circ g$. Since $\tau\iota=f$ (the type of $W\cap\Sph$ is $M_W$), $\tau gf\simeq\mu f\simeq\tau\iota$. As $\Gr_n^d$ is a CW complex and $\tau$ a weak equivalence, $gf\simeq\iota$. So $rf=\rho gf\simeq\rho\iota=\mathrm{id}$.
\end{proof}

\begin{theoremM}\label{thm:M}
Let $n\ge4$. Hypotheses \ref{H:basin} and (HM) imply that the realization map $f:\Gr_n^d\to\|\AM_n^d\|$ is a homotopy equivalence. In particular they imply Conjecture~I.
\end{theoremM}
\begin{proof}
Write $X=\|\AM_n^d\|$, $G=\Gr_n^d$ and $D=d(n-1-d)$, and let $r$ be as in Corollary~\ref{cor:Bprime}.

\emph{Double covers.} Since $n-1\ge3$, the oriented Grassmannian $\widetilde G\to G$ is the universal cover; it is a connected double cover, and $\widetilde G$ is a closed orientable manifold. Let $\widehat X=r^*\widetilde G$. Then $f^*\widehat X\cong(rf)^*\widetilde G\cong\widetilde G$ is connected, so $\widehat X$ is connected. Let $\hat f:\widetilde G\to\widehat X$ and $\hat r:\widehat X\to\widetilde G$ be lifts of $f$ and $r$. Lifting a homotopy $rf\simeq\mathrm{id}$ shows that $\hat r\hat f$ is homotopic to a deck transformation; composing $\hat f$ with it, $\hat r\hat f\simeq\mathrm{id}$. As a finite cover of a compact $\Z$-homology manifold, $\widehat X$ is a closed connected $\Z$-homology manifold of dimension $D$.

\emph{Degree.} $H_D(\widetilde G;\Z)\cong\Z$ and $\hat r_*\hat f_*=\mathrm{id}$, so $H_D(\widehat X;\Z)\ne0$. Hence $\widehat X$ is orientable \cite[Ch.~8]{Munkres}, and $\deg\hat r\cdot\deg\hat f=1$, i.e.\ $\deg\hat f=\pm1$.

\emph{Homology.} By Poincar\'e duality \cite[Ch.~8]{Munkres}, $\hat f_!=PD\circ\hat f^*\circ PD^{-1}$ satisfies $\hat f_*\hat f_!=\pm\mathrm{id}$. Indeed, for $\alpha=\alpha'\cap[\widehat X]$ we have $\hat f_*(\hat f^*\alpha'\cap[\widetilde G])=\alpha'\cap\hat f_*[\widetilde G]=\pm\alpha$. So $\hat f_*$ is surjective in integral homology. It is injective because $\hat r_*\hat f_*=\mathrm{id}$.

\emph{Fundamental group.} A map of degree $\pm1$ is surjective on $\pi_1$. Otherwise $\hat f$ lifts to the cover $X'\to\widehat X$ associated with its image. If the index $m\ge2$ is finite, $\deg\hat f=m\deg\hat f'$, which is impossible. If the index is infinite, $X'$ is not compact, $H_D(X';\Z)=0$, and $\hat f_*[\widetilde G]=0$. Since $\pi_1(\widetilde G)=1$, $\widehat X$ is simply connected.

\emph{Conclusion.} $\hat f$ is an integral homology isomorphism between simply connected CW complexes, hence a homotopy equivalence. Both covers are universal double covers and $\hat f$ is equivariant, so $f$ induces isomorphisms $\pi_1(G)=\Z/2\to\pi_1(X)=\Z/2$ and $\pi_i(G)\to\pi_i(X)$ for $i\ge2$. Whitehead's theorem concludes.
\end{proof}

\noindent Theorem~M needs neither \ref{H:type} nor \ref{H:degree}: the flow only has to converge, not to preserve or follow combinatorial types, so the obstruction of Theorem~O below does not apply to it. The work is shared as follows.
\begin{itemize}
\item The flow supplies the integral and $\pi_1$ information, through the retraction $r$.
\item Poincar\'e duality supplies the ``upper half''.
\item Without the flow, Theorem~W and Corollary~\ref{cor:manifold} already give the mod~$2$ isomorphism.
\end{itemize}
(HM) is local and combinatorial: every open interval, and every lower and upper interval, of $\AM_n^d$ must be a $\Z$-homology sphere of the right dimension. It holds over $\F_2$ and $\F_3$ for $n\le5$, and for $d\in\{1,n-2\}$ by Theorem~A.

\subsection{The degree route}
Let $U^\chi_M=\{\Sph'\in\fat_n^d:\tau(\Sph')\ge_\chi M\}$.

\begin{proposition}\label{prop:Schi}
Assume (W). Then each $U^\chi_M$ is open and contractible, so $\tau:\fat_n^d\to\AM_n^{d,\chi}$ is a weak equivalence. Moreover $\overline\Phi:\fat_n^d\to\Pl/\pm$ is a weak equivalence, and $\overline\Phi\circ\iota\simeq p$.
\end{proposition}
\begin{proof}
Openness follows from Lemma~\ref{lem:degclosure}, as in Lemma~\ref{lem:tau-cont}. By (W), $U^\chi_M\subseteq U_M\subseteq R(M)$, and the contraction of Theorem~S preserves types for $t<1$, so it contracts $U^\chi_M$ within itself. McCord's theorem applies.

For $q$ in the cell of $\pm\chi^*_M$, the stars $\mathrm{St}(q)$ form a basis-like cover of $\Pl/\pm$ by contractible open sets (Proposition~\ref{prop:plucker}), and by Theorem~P(b) the preimage of the star of $M$ is $U^\chi_M$. Again McCord's theorem applies.

For $W\in\Gr_n^d$, $\hat\pi_B|_{W\cap\Sph}$ is the radial projection of the linear map $\pi_B|_W$. So $\Phi(W\cap\Sph)$ and $p(W)$ lie in the same convex open cell, and the straight-line homotopy descends to $\Pl/\pm$.
\end{proof}

\begin{theorem}[Theorem E]\label{thm:E}
Assume \ref{H:basin}, \ref{H:degree} and (W). Then $\iota:\Gr_n^d\hookrightarrow\fat_n^d$ is a weak homotopy equivalence. Consequently $\|\AM_n^d\|\simeq\|\AM_n^{d,\chi}\|\simeq G(d,n-1)$, and $p:\Gr_n^d\to\Pl/\pm$ is a weak equivalence.
\end{theorem}
\begin{proof}
Let $\rho=\Lambda\sigma$. Homotope $\overline\Phi$ to $p\circ\rho$, as maps $\fat_n^d\to\Pl/\pm$, in three stages:
\begin{enumerate}
\item from $\overline\Phi(\Sph')$ to $\overline\Phi(\sigma\Sph')$, by a segment in a convex open cell (\ref{H:degree}(a));
\item along $t\mapsto\overline\Phi(\Sph'_t)$, $t\in[0,\infty]$. This path lies in $\Pl/\pm$ by \ref{H:degree}(c) and Theorem~P, and it is continuous by \ref{H:degree}(b), since degrees are locally constant away from $K_B$;
\item from $\overline\Phi(\rho(\Sph'))$ to $p(\rho(\Sph'))$, as in Proposition~\ref{prop:Schi}.
\end{enumerate}
Since $\overline\Phi$ is a weak equivalence, $\rho_*$ is injective on homotopy groups. Since $\rho\iota=\mathrm{id}$, $\rho_*$ is surjective. So $\iota_*=\rho_*^{-1}$ is an isomorphism. Theorem~S and Proposition~\ref{prop:Schi} conclude.
\end{proof}

\begin{theorem}[Theorem E$^\flat$]\label{thm:Eflat}
Assume \ref{H:basin}, (H4$^\flat$) and (W). Then $\iota:\Gr_n^d\hookrightarrow\fat_n^d$ is a weak homotopy equivalence, and Conjecture~I holds.
\end{theorem}
\begin{proof}
Let $\widehat S$ be the finite space of nonzero sign patterns on $\binom{[n]}\ell$ modulo $\pm$, $P_\chi=\{\pm\chi^*_M\}/\pm\subseteq\widehat S$, $Y=\fat_n^d\times[0,\infty]$ and $g(\Sph',t)=\psi_t(\Sph')$.
\begin{itemize}
\item \emph{$g$ is continuous.} Near a point of $Y$ the nonzero entries of the pattern are locally constant by \ref{H:degree}(b), and zero entries may take any value nearby, which is what the minimal open sets allow. At $t=0$ use \ref{H:degree}(a) and Lemma~\ref{lem:degclosure}; at $t=\infty$, use convergence to $\rho(\Sph')$.
\item \emph{Fibre lemma.} Let $j:P_\chi\hookrightarrow P_\chi\cup g(Y)$. Its fibres over the down-sets of $s$ are the posets $F(s)$. They are contractible for $s\in P_\chi$, since they have a maximum, and for $s\in g(Y)$ by (H4$^\flat$). By Lemma~\ref{lem:fibre} (applied to the opposite posets) and the naturality of $\mu$, $j$ is a weak equivalence.
\item \emph{Cylinder.} $g\circ i_0=j\circ\tau$ by Theorem~P, and $g\circ i_\infty=j\circ\chi\circ\rho$ with $\chi(W)=\pm\chi^*_{M_W}$. So $\tau_*=\chi_*\rho_*$ up to base points. Since $\tau$ is a weak equivalence (Proposition~\ref{prop:Schi}), $\rho_*$ is injective, and as before $\iota_*$ is an isomorphism.
\end{itemize}
\end{proof}

\noindent (H4$^\flat$) is a genuine condition. For $(5,2)$ there are $1600$ sign patterns whose fibre is nonempty and disconnected. For the acyclic cases $(4,1)$ and $(5,1)$, every fibre is empty or contractible (computer).

\subsection{The obstruction}
Let $\mathcal D\subseteq\Gr_n^d$ be the closed set of $W$ at which the Pl\"ucker coordinates $p_B$ vanishing at $W$ do not have linearly independent differentials on $\Gr_n^d$. Examples are configurations with two coincident points and, for $(8,3)$, with five coplanar points.

\begin{proposition}[start and end]\label{prop:startend}
\begin{enumerate}[label=(\alph*)]
\item Assume \ref{H:degree}(a),(b). If every completion of $\chi^*_{\tau(\Sph')}$ (e.g.\ if $\tau(\Sph')$ is uniform, or $\chi^*_{\tau(\Sph')}$ has at most one zero) is $\pm\chi^*_{M'}$ for some $M'\in\AM_n^d$, then \ref{H:degree}(c) holds along $\Sph'_t$ for $t\in[0,t_0)$, for some $t_0>0$.
\item Assume \ref{H:basin} and \ref{H:degree}(b). If $\rho(\Sph')\notin\mathcal D$, there is $T$ such that $\psi_t(\Sph')$ is the dual chirotope of a realizable oriented matroid for $t\ge T$.
\end{enumerate}
\end{proposition}
\begin{proof}
(a) Let $\delta=\min\{\mathrm{dist}(\sigma\Sph',K_B):\chi^*_{\tau(\Sph')}(B)\ne0\}>0$. For $t$ small $\Sph'_t$ stays within $\delta/2$ of $\sigma\Sph'$, so these degrees do not change, and the pattern is a completion (Lemma~\ref{lem:completion} for one zero).

(b) Let $W=\rho(\Sph')$. For $B$ with $p_B(W)\ne0$, $C^0$ convergence gives $\deg_B=\sgn p_B(W)$ for large $t$. The other entries form a completion of $\chi_W$. Since $W\notin\mathcal D$, the map $W'\mapsto(p_B(W'))_{p_B(W)=0}$ is a submersion at $W$, hence open, so every completion is realized by some $W'$ near $W$.
\end{proof}

\begin{proposition}[curves]\label{prop:curves}
Let $k=1$ and assume \ref{H:basin}. If each coordinate $x_e$ has at most two simple zeros on $\sigma(\Sph')$, then \ref{H:degree}(c) holds for all $t$.
\end{proposition}
\begin{proof}
By Lemma~\ref{lem:coordinates} and Angenent's theorem \cite{Ang88}, the number of zeros of $x_e$ does not increase and drops at multiple zeros. An odd function has at least two zeros, so $x_e$ keeps one antipodal pair of simple zeros. Then $\psi_t$ is the chirotope of the cyclic arrangement of these pairs, a rank $2$ configuration; collisions produce parallel elements, which still give a chirotope. It is totally cyclic, because otherwise all $x_e$ would be positive at some point, contradicting $\sum_ex_e=0$.
\end{proof}

\begin{theoremO}\label{thm:O}
Assume \ref{H:basin}, \ref{H:degree} and, for a uniform $M\in\AM_n^d$, the genericity condition
\begin{itemize}
\item[(G)] some $\Sph'\in\tau^{-1}(M)$ has a trajectory with finitely many crossing times, each simple (one $B$, one antipodal pair, transversal), and $\rho(\Sph')\notin\mathcal D$.
\end{itemize}
Then $\chi_M$ is joined to the chirotope of a realizable oriented matroid by a sequence of mutations through uniform chirotopes. If (G) holds for every uniform $M\in\AM_n^d$, the mutation graph of uniform oriented matroids of rank $d+1$ on $[n]$ is connected.
\end{theoremO}
\begin{proof}
Between crossing times the pattern is constant and uniform. At a simple crossing at $B$, $\deg_B$ changes by $\pm2$. The pattern passes from $\psi^-$ through $\psi^0$, which has a single zero at $B$, to $\psi^+$, which is either $\psi^-$ (when $|\deg_B|>1$) or $\psi^-$ with $B$ reversed. By \ref{H:degree}(c) both are chirotopes, so the second case is a mutation. The pattern starts at $\pm\chi^*_M$ (Theorem~P and \ref{H:degree}(a)) and ends at a realizable uniform chirotope (Proposition~\ref{prop:startend}(b)). Duality preserves mutations, since $\chi(E\setminus B)=\pm\chi^*(B)$ with a sign depending only on $B$.

For connectivity, note the following.
\begin{itemize}
\item Every uniform oriented matroid has an acyclic reorientation. Reorientation is an automorphism of the mutation graph that preserves realizability.
\item Uniform realizable chirotopes form a single mutation class \cite{RS88}: a generic path in the connected manifold $\Gr_n^d$ crosses the hypersurfaces $\{p_B=0\}$ one at a time and transversally.
\end{itemize}
\end{proof}

\begin{corollary}\label{cor:GK}
For $d\ge6$ and $n\ge d+18$, Hypotheses \ref{H:basin}, \ref{H:degree} and (G) for all uniform $M$ cannot hold simultaneously.
\end{corollary}
\begin{proof}
By Gu and Knauer \cite[Thm.~1.6]{GK26}, for every $r\ge7$ and $n\ge r+17$ the mutation graph of uniform oriented matroids of rank $r$ on $[n]$ is disconnected. Take $r=d+1$ in Theorem~O.
\end{proof}

\begin{remark}[how the obstruction enters]\label{rem:obstruction}
The mechanism is transparent. A generic trajectory in Pl\"ucker space crosses only codimension-one walls of the orthant decomposition. If two uniform chirotopes lie in different mutation classes, any path between them through valid cells must pass through a cell of codimension at least $2$, i.e.\ a degenerate oriented matroid with at least two vanishing bases. A generic trajectory therefore has to pass through \emph{invalid} patterns.

Corollary~\ref{cor:GK} rules out every argument that follows the type of the sphere generically along the flow. The numerical observation that types change only by mutations (Appendix~\ref{app:numerics}) cannot reflect the general behaviour. It also explains why this paper turns to spherical separoids. The arguments that survive are of three kinds:
\begin{itemize}
\item those using only the limit of the flow (Theorem~M, Corollary~\ref{cor:Bprime});
\item those that tolerate invalid patterns (Theorem~\ref{thm:Eflat});
\item those that replace oriented matroids by spherical separoids (Section~\ref{subsec:flow-separoids}).
\end{itemize}
Note that (H4$^\flat$) does not imply connectivity: an invalid pattern whose fibre consists of degenerate chirotopes can connect different mutation classes.
\end{remark}

\subsection{Spherical separoids along the flow}\label{subsec:flow-separoids}
For $y=(\Sph',t)\in Y=\fat_n^d\times[0,\infty]$, let $\bar\Sp(y)$ be the separoid whose Radon partitions are the signed sets lying above the sign vectors of points of $\Sph''_{t'}$ for \emph{all} $(\Sph'',t')$ in arbitrarily small neighbourhoods of $y$ (a $\limsup$). Here $\Sph''_\infty=\rho(\Sph'')$. Since there are finitely many signed sets, some neighbourhood $U$ of $y$ realizes the intersection, and $\bar\Sp(y')\ge\bar\Sp(y)$ for $y'\in U$. So $\bar\Sp$ is continuous for the Alexandrov topology. Let $\mathcal B^{\F}$ be the poset of all $\F$-spherical acyclic separoids on $[n]$ with sphere dimension $\ell-1$, with no dimension condition. It is admissible when $\F=\F_2$.

\begin{proposition}\label{prop:factor}
Assume \ref{H:basin}, \ref{H:degree}(b), (a$^{\Sep}$) and (H4$^{\Sep}$). Then the inclusion $j:\|\AM_n^d\|\hookrightarrow\|\mathcal B^{\F}\|$ is homotopic to $f\circ r$, where $f$ is the realization map and $r$ is as in Corollary~\ref{cor:Bprime}. Consequently:
\begin{enumerate}[label=(\alph*)]
\item $j\simeq j\circ e$, where $e=f_{\AM}\circ r$ is a homotopy idempotent of $\|\AM_n^d\|$ with image $\Gr_n^d$. Inside the space of spherical separoids, the MacPhersonian collapses onto the Grassmannian.
\item $j_*$ kills the complement $\ker r_*$ of $f_*H_*(\Gr_n^d)$ in $H_*(\|\AM_n^d\|)$. For $\F=\F_2$ coefficients, $\ker j_*=\ker r_*$ (by Theorem~W).
\item If moreover (D) holds for $\mathcal B^{\F}$, then $\|\AM_n^d\|\simeq G(d,n-1)$.
\end{enumerate}
\end{proposition}
\begin{proof}
$\bar\Sp:Y\to\mathcal B^{\F}$ is continuous by the above and lands in $\mathcal B^{\F}$ by (H4$^{\Sep}$).
\begin{itemize}
\item \emph{$t=0$.} $\bar\Sp(\Sph',0)=\tau(\Sph')$: the Radon set contains $\Rad(\tau(\Sph'))$ by (a$^{\Sep}$), and is contained in it by semicontinuity (Lemma~\ref{lem:semicont}) and \ref{H:degree}(b).
\item \emph{$t=\infty$.} Likewise $\bar\Sp(\Sph',\infty)=M_{\rho(\Sph')}$.
\end{itemize}
So $j\circ\tau\simeq f\circ\rho$ as maps $\fat_n^d\to\mathcal B^{\F}$, the homotopy being the cylinder $Y$. Composing with $g$ as in Corollary~\ref{cor:Bprime}, and identifying maps from CW complexes into finite spaces with their lifts through $\mu$ (Section~\ref{sec:tools}), gives $j\simeq j\tau g\simeq f\rho g=fr$. Then (a) follows from $f=jf_{\AM}$. For (b), $j_*=f_*r_*$, and $f_*$ is injective over $\F_2$ by Theorem~W. For (c), if $j$ is a weak equivalence then $f_{\AM}r\simeq\mathrm{id}$, and together with $rf_{\AM}\simeq\mathrm{id}$ this shows that $f_{\AM}$ is a homotopy equivalence.
\end{proof}

\noindent Proposition~\ref{prop:factor} is all that the flow can say about $\Sep_n^d$ when it starts in $\fat_n^d$. Separoids not reached by spheres of oriented matroid type, such as the $480$ separoids of $\Sep_5^2$ incomparable with every oriented matroid, are invisible to it. Conversely, the flow cannot contract the spaces $\fat(\Sp)$ of Corollary~\ref{cor:D-spheres} through $\Gr_n^d$, because $\Omega(M)=\emptyset$ for non-realizable uniform $M$. The flow and (D) are therefore complementary. Under \ref{H:basin} and (HM), Theorem~M and Corollary~\ref{cor:DF2} give
\[\|\Sep_n^d\|\simeq G(d,n-1)\iff\text{(D) for }\Sep_n^d,\qquad\text{and, mod }2,\ \iff\ \textstyle\sum_ib_i(\|\Sep_n^d\|;\F_2)=\binom{n-1}{d}.\]

\subsection{Summary of the proved implications}
\begin{center}\footnotesize
\begin{tabular}{@{}p{0.34\textwidth}p{0.44\textwidth}l@{}}
\toprule
hypotheses & conclusion & reference\\
\midrule
none & (H1), i.e.\ $\fat_n^d\simeq_w\|\AM_n^d\|$ & Thm.~S\\
none & $H^*(G;\F_2)$ is a ring retract of $H^*(\|\mathcal B\|;\F_2)$ & Thm.~W\\
none & Conjectures I and II for $d=n-2$; Conjecture I for $d=1$ & Thms.~A, F\\
\ref{H:type} & Conjecture I & Thm.~\ref{thm:B}\\
\ref{H:basin} & $G$ homotopy retract of $\|\AM_n^d\|$ & Cor.~\ref{cor:Bprime}\\
\ref{H:basin}, (HM) & Conjecture I & Thm.~M\\
\ref{H:basin}, \ref{H:degree}, (W) & Conjecture I; $\Pl/\pm\simeq G$ & Thm.~\ref{thm:E}\\
\ref{H:basin}, (H4$^\flat$), (W) & Conjecture I & Thm.~\ref{thm:Eflat}\\
\ref{H:basin}, \ref{H:degree}, (G) & uniform mutation graph connected & Thm.~O\\
\ref{H:basin}, \ref{H:degree}, (G), \cite{GK26} & contradiction for $d\ge6$, $n\ge d+18$ & Cor.~\ref{cor:GK}\\
\ref{H:basin}, \ref{H:degree}(b), (a$^{\Sep}$), (H4$^{\Sep}$) & $j\simeq f\circ r$ & Prop.~\ref{prop:factor}\\
\quad $+$ (D) & Conjecture I & Prop.~\ref{prop:factor}(c)\\
(HM)$_{\F_2}$ & $H^*(\|\AM_n^d\|;\F_2)\cong H^*(G;\F_2)$ & Cor.~\ref{cor:manifold}\\
weakly contractible $\fat(\Sp)$ for all $\Sp$ & (D) for $\Sep_n^d$ & Cor.~\ref{cor:D-spheres}\\
\bottomrule
\end{tabular}
\end{center}
% =====================================================================
\part*{IV. Open problems}
\section{Open problems and conjectures}\label{sec:open}

\subsection{The two conjectures}
Conjecture~I is open in the range $3\le d\le n-5$. Conjecture~II splits into two statements:
\begin{enumerate}[label=(II.\alph*)]
\item the realization map $G(d,n-1)\to\|\Sep_n^d\|$ is a homotopy equivalence;
\item (D) the inclusion $\AM_n^d\hookrightarrow\Sep_n^d$ is a homotopy equivalence.
\end{enumerate}
Theorem~W proves the ``lower half'' of (II.a) in mod~$2$ cohomology. Theorem~F proves both for $d=n-2$. The computations prove both for $n\le5$ for the subposet of spherical Radon separoids.

\begin{conjecture}[upper bound]\label{conj:betti}
$\sum_ib_i(\|\Sep_n^d\|;\F_2)\le\binom{n-1}{d}$.
\end{conjecture}

\noindent By Corollary~\ref{cor:betti}, Conjecture~\ref{conj:betti} is equivalent to the mod~$2$ form of (II.a), namely $f^*$ being an isomorphism. A natural route is an acyclic matching on $\Sep_n^d$ (discrete Morse theory) with exactly $\binom{n-1}{d}$ critical cells, one for each Schubert cell. Proposition~\ref{prop:nolocal} shows that no proof can proceed by contracting preimages of open stars under the realization map.

\begin{conjecture}[spheres in Radon complexes]\label{conj:fibres}
For every $\Sp\in\Sep_n^d$, the space $\fat(\Sp)$ of piecewise geodesic spheres of oriented matroid type contained in the Radon complex $R(\Sp)$ is weakly contractible.
\end{conjecture}

\noindent By Corollary~\ref{cor:D-spheres}, Conjecture~\ref{conj:fibres} implies (II.b). It is a statement about spaces of embedded spheres inside an open set homotopy equivalent to a sphere, with the extra rigidity of the cell structure. For $\Sp\in\AM_n^d$ it is Theorem~S.

\begin{conjecture}[local MacPherson conjecture]\label{conj:HM}
$\|\AM_n^d\|$ is a closed $\Z$-homology manifold of dimension $d(n-1-d)$. Equivalently, $\AM_n^d$ is graded of rank $d(n-1-d)$, and all its lower, upper and open intervals are $\Z$-homology spheres of the appropriate dimensions.
\end{conjecture}

\noindent This holds over $\F_2$ and $\F_3$ for $n\le5$. By Corollary~\ref{cor:manifold}, its $\F_2$-form implies $H^*(\|\AM_n^d\|;\F_2)\cong H^*(G(d,n-1);\F_2)$. By Theorem~M, together with \ref{H:basin}, it implies Conjecture~I. It can be tested locally in the first open case. For instance, for the certified non-realizable $\chi_1$ of Proposition~\ref{prop:nonreal}, a maximal element of $\AM_8^3$, it predicts that $\|(\AM_8^3)_{<\chi_1}\|$ is a homology $11$-sphere. This is a finite computation.

\subsection{The flow}
\begin{conjecture}[basin]\label{conj:basin}
Hypothesis~\ref{H:basin} holds: after a suitable smoothing, the mean curvature flow of every odd piecewise geodesic sphere of oriented matroid type exists for all time and converges to a great sphere, continuously in the initial datum.
\end{conjecture}

\noindent Near $\Gr_n^d$ this is a small-curvature convergence statement (Section~\ref{sec:conditional}). The geodesic condition on the cells is essential here: a knotted sphere of the topological model (Remark~\ref{rem:S}) cannot converge to a great sphere through embedded spheres. A proof for $(n,d)=(8,3)$, where the spheres are $3$-spheres in $S^6$, would settle Conjecture~I there, given Conjecture~\ref{conj:HM}.

\begin{conjecture}[spherical shadow]\label{conj:H4S}
(H4$^{\Sep}$) holds: the separoid shadow of the flow is spherical at all times.
\end{conjecture}

\noindent Corollary~\ref{cor:GK} shows that the analogous statement for oriented matroids is false. Numerically (Appendix~\ref{app:numerics}), the flow leaves the oriented matroids and even the Radon separoids, but stays spherical in all $11$ further runs that we checked.

\begin{question}
Along a generic trajectory, is every invalid degree pattern $\psi$ visited by the flow such that $F(\psi)$ is contractible? This is (H4$^\flat$). For the invalid pattern of Appendix~\ref{app:numerics}, $F(\psi)$ is nonempty, but its contractibility is unknown.
\end{question}

\subsection{Combinatorial questions}
\begin{question}[$d=1$]\label{q:d1}
Let $\Sp\in\Sep_n^{1,\F}$. By Proposition~\ref{prop:alexander} and the argument of Theorem~F, $L(\Sp)$ is the disjoint union of two $\F$-acyclic halves $D\sqcup-D$. An oriented matroid in $\AM_n^1$ is a nonempty chain $c$ of nonempty proper subsets of $[n]$, up to complementation of all its members. Unwinding the definitions, $M_c\ge\Sp$ iff, for one of the two orientations of $c$, $c$ meets every Boolean interval
\[J_Z=[Z^-,\,[n]\setminus Z^+],\qquad Z\in D.\]
Is the complex of chains of the Boolean lattice meeting all the intervals $J_Z$, $Z\in D$, always $\F$-acyclic? By Lemma~\ref{lem:fibre}(b), a positive answer gives (II.b) for $d=1$ in $\F$-homology. For oriented matroids the intervals are singletons, and the complex is a closed star.
\end{question}

\begin{question}[degrees of Radon complexes]
For $\Sp\le\T$ in $\Sep_n^d$, is the inclusion $R(\T)\subseteq R(\Sp)$ a homotopy equivalence? Lemma~\ref{lem:gysin}(c) shows that it has odd degree. A positive answer would make the tautological fibration an honest spherical quasifibration over $\|\Sep_n^d\|$, with an integral Euler class and rational Pontryagin classes. That would give integral and rational analogues of Theorem~W.
\end{question}

\begin{question}[orders]
Do the separoid order, the vector order and the weak-map order coincide on $\AM_n^d$? The first two coincide for $n\le5$. Our results use only the separoid order, except Theorems~\ref{thm:E} and~\ref{thm:Eflat}, which use (W).
\end{question}

\begin{question}[Whitney duality]
The separation complexes $L(\Sp)$ form a second tautological family of mod~$2$ $(d-1)$-spheres, covariant in $\Sp$. Is $w(R)\,w(L)=1$ in $H^*(\|\mathcal B\|;\F_2)$, as for $\gamma\oplus\gamma^\perp=\one^\perp$?
\end{question}

\begin{question}[computations]
Extend Theorem~\ref{thm:C} to $n=6,7$. Decide (II.b) and Conjecture~\ref{conj:betti} for all spherical separoids (without the Radon condition) for $n=5$.
\end{question}
% =====================================================================
\appendix
\section{Numerical evidence}\label{app:numerics}

The experiments in this appendix are numerical: one discretization per run, checks at discrete times, finitely many initial data. They are evidence, not proof, and nothing in the body of the paper depends on them. All of them concern $(n,d)=(8,3)$, the first open case, where the spheres are $3$-spheres in $\Sph\cong S^6$ and $\ell=4$, except Observation~\ref{obs:curves}.

\subsection{Method}
For a type $M\in\AM_8^3$, we evolved odd piecewise geodesic spheres of type $M$ by mean curvature flow, after triangulating them. The initial spheres were either the barycentric sphere $\Sph_M$ ($928$ vertices, $5376$ tetrahedra) or spheres obtained from it by moving every vertex to a random point of its cell (Lemma~\ref{lem:moving}). We used the following discretization:
\begin{itemize}
\item a semi-implicit Dziuk scheme $(\mathbf M+\Delta t\,\mathbf K)X^{\mathrm{new}}=\mathbf MX^{\mathrm{old}}$ with the $P_1$ mass and stiffness matrices, followed by projection to $\Sph$;
\item an optional tangential redistribution;
\item exact antipodal symmetrization.
\end{itemize}
Along the flow we recorded three things:
\begin{itemize}
\item the \emph{circuit census}, i.e.\ the number of intersection points of the sphere with each great $3$-sphere $L_S=\Sph\cap\R^S$, $|S|=5$ (exactly one antipodal pair for every $S$ for a sphere of uniform oriented matroid type);
\item the \emph{degree pattern} $(\sgn\deg_B)_{|B|=4}$, computed as signed counts of preimages of random directions;
\item the \emph{separoid} of the sphere, with its sphericity tested over $\F_2$ and $\F_3$ via Proposition~\ref{prop:alexander}.
\end{itemize}

\subsection{Observations}
\begin{observation}[degree theorem]\label{obs:degree}
For $\chi_1,\chi_2$ (Proposition~\ref{prop:nonreal}) and two realizable controls, on $\Sph_M$ and on fourteen spheres with randomly moved vertices, all $70$ degrees satisfied $|\deg_B|=1$ and $\deg_B=\varepsilon\chi^*_M(B)$, as predicted by Theorem~P. The projections are far from injective: a generic point has up to $67$ preimages.
\end{observation}

\begin{observation}[convergence]\label{obs:n8}
From $\Sph_M$, for $M=\chi_1,\chi_2$ and both controls:
\begin{itemize}
\item the flow converged to a great $3$-sphere, with exponential rate $\approx12$, the first positive eigenvalue $l(l+k-1)-k$ ($l=k=3$) of the Jacobi operator on odd normal fields;
\item at every check the sphere had uniform oriented matroid type;
\item the type changed only by mutations, $8$ and $5$ of them for $\chi_1$ and $\chi_2$, and $3$ for each control. Starting from $\chi_1,\chi_2$, the first mutation already produced a realizable oriented matroid.
\end{itemize}
Without symmetrization, a mode growing at rate $\approx3=k$ appears, which is the direction of the round point.
\end{observation}

\begin{observation}[violation of \ref{H:degree}(c)]\label{obs:violation}
For one sphere of type $\chi_1$ with randomly moved vertices, the following sequence occurs along the flow.
\begin{enumerate}[label=(\alph*)]
\item The census breaks: some $L_S$ meets the sphere in $6$ points. This is the saddle mechanism: $u=x_1^2-3x_2^2-4t$ solves $\partial_tu=\Delta u$, and $u(x_1,0,t)$ acquires two zeros at $t=0$.
\item The sphere then crosses $K_B$ for a basis $B$ that is \emph{not} mutable in the current chirotope, and the degree pattern violates a $3$-term Grassmann--Pl\"ucker relation.
\item Shortly afterwards it crosses a second $K_{B'}$. The pattern becomes valid again, and later the type is restored.
\end{enumerate}
The window of invalidity was found in four discretizations, with the same pair of crossings in all of them.
\begin{center}\small
\begin{tabular}{@{}lc@{}}
\toprule
scheme & invalid window\\
\midrule
$928$ vertices, $\Delta t=2\cdot10^{-3}$, redistribution & $[0.028,0.036)$\\
$928$ vertices, $\Delta t=5\cdot10^{-4}$, redistribution & $[0.010,0.014)$\\
refined ($7232$ vertices), $\Delta t=10^{-3}$, redistribution & $[0.026,0.038)$\\
refined, $\Delta t=10^{-3}$, plain Dziuk scheme & $[0.070,0.077)$\\
\bottomrule
\end{tabular}
\end{center}
In the window, the separoid of the sphere is not a Radon separoid (a $5$-set carries two or three minimal Radon partitions), but it is spherical over $\F_2$ and $\F_3$. For the invalid pattern $\psi$, the fibre $F(\psi)$ is nonempty: it contains dual chirotopes of degenerate realizable configurations.
\end{observation}

\begin{observation}[spherical shadows]\label{obs:separoid}
In $11$ further flows from random piecewise geodesic spheres of random uniform types and of types $\chi_1,\chi_2$, the sphere left oriented matroid type at $22$ checks. At $2$ of them its separoid was not a Radon separoid, and at all $22$ it was spherical. One of these flows produced a second violation of \ref{H:degree}(c), again with a non-Radon spherical separoid.
\end{observation}

\begin{observation}[random curves]\label{obs:curves}
Sphericity is not automatic for separoids of odd spheres. Among random odd trigonometric curves of degree $3$ to $9$ in $\Sph\subset\R^5$ ($(n,d)=(5,2)$, $400$ curves per degree), between $6\%$ and $31\%$ have a non-spherical separoid. Sphericity and validity of the degree pattern are independent: for degree $9$, $76$ curves have a spherical separoid and an invalid pattern, and $65$ the opposite. The curve of Example~\ref{ex:badpattern} has a non-Radon spherical separoid.
\end{observation}

\noindent Taken together, these observations support Conjecture~\ref{conj:H4S}, and they are consistent with Corollary~\ref{cor:GK}. Typical trajectories cross walls at mutations, but some trajectories pass through invalid patterns, and when they do, the separoid stays spherical.
% =====================================================================
\section*{Acknowledgments}
This work was developed with the assistance of Claude, an AI model by Anthropic. It contributed to the development and checking of the arguments, wrote the code and carried out the computations of Section~\ref{sec:computations} and Appendix~\ref{app:numerics}. The author is responsible for the content.

\end{document}